\documentclass[11pt]{article}
\usepackage{latexsym,amssymb}
\usepackage{psfrag}
\usepackage{amsmath}
\usepackage{color}
\usepackage{subfig}
\usepackage[pdftex]{graphicx}

\usepackage{tikz}
\usetikzlibrary{intersections, calc, angles, arrows.meta}
\usetikzlibrary{patterns}

\newtheorem{Theorem}{\bf Theorem}[section]
\newtheorem{Lemma}{\bf Lemma}[section]
\newtheorem{Proposition}{\bf Proposition}[section]
\newtheorem{Corollary}{\bf Corollary}[section]
\newtheorem{Remark}{\bf Remark}[section]
\newtheorem{Example}{\bf Example}[section]
\newtheorem{Definition}{\bf Definition}[section]

\newenvironment{theorem}{\begin{Theorem}$\!\!\!$}{\end{Theorem}}
\newenvironment{lemma}{\begin{Lemma}$\!\!\!$}{\end{Lemma}}

\newenvironment{corollary}{\begin{Corollary}$\!\!\!$}{\end{Corollary}}
\newenvironment{remark}{\begin{Remark}$\!\!\!$}{\end{Remark}}

\newenvironment{definition}{\begin{Definition}$\!\!\!$}{\end{Definition}}

\def\XXint#1#2#3{{\setbox0=\hbox{$#1{#2#3}{\int}$}
\vcenter{\hbox{$#2#3$}}\kern-.5\wd0}}

\numberwithin{equation}{section}

\allowdisplaybreaks[4]

\begin{document}

\title{Initial traces and solvability of Cauchy problem\\
to a semilinear parabolic system}
\author{
Yohei Fujishima and Kazuhiro Ishige
}
\date{}
\maketitle
\begin{abstract}
  Let $(u,v)$ be a solution to a semilinear parabolic system
  $$
  \mbox{(P)}
  \qquad
  \left\{
  \begin{array}{ll}
  \partial_t u=D_1\Delta u+v^p\quad & \quad\mbox{in}\quad{\bf R}^N\times(0,T),\vspace{3pt}\\
  \partial_t v=D_2\Delta v+u^q\quad &  \quad\mbox{in}\quad{\bf R}^N\times(0,T),\vspace{3pt}\\
  u,v\ge 0 &  \quad\mbox{in}\quad{\bf R}^N\times(0,T),\vspace{3pt}\\
  (u(\cdot,0),v(\cdot,0))=(\mu,\nu) & \quad\mbox{in}\quad{\bf R}^N,
  \end{array}
  \right.
  $$
  where $N\ge 1$, $T>0$, $D_1>0$, $D_2>0$, $0<p\le q$ with $pq>1$ and $(\mu,\nu)$ is a pair of
  Radon measures or nonnegative measurable functions in ${\bf R}^N$.
  In this paper we study qualitative properties of the initial trace of the solution~$(u,v)$
  and obtain necessary conditions on the initial data $(\mu,\nu)$ for the existence of solutions to problem~(P).
\end{abstract}
\vspace{50pt}
\noindent Addresses:

\smallskip
\noindent Y.~F.:  Department of Mathematical and Systems Engineering, Faculty of Engineering, Shizuoka University,
3-5-1, Johoku, Hamamatsu, 432-8561, Japan.\\
\noindent
E-mail: {\tt fujishima@shizuoka.ac.jp}\\

\smallskip
\noindent K.~I.:  Graduate School of Mathematical Sciences, The University of Tokyo,
3-8-1 Komaba, Meguro-ku, Tokyo 153-8914, Japan.\\
\noindent
E-mail: {\tt ishige@ms.u-tokyo.ac.jp}
\newpage

\section{Introduction}
Consider a semilinear parabolic system
\begin{equation}
\label{eq:1.1}
\left\{
\begin{array}{ll}
\partial_t u=D_1\Delta u+v^p & \quad\mbox{in}\quad{\bf R}^N\times(0,T),\vspace{3pt}\\
\partial_t v=D_2\Delta v+u^q & \quad\mbox{in}\quad{\bf R}^N\times(0,T),\vspace{3pt}\\
u,v\ge 0 & \quad\mbox{in}\quad{\bf R}^N\times(0,T),
\end{array}
\right.
\end{equation}
where $N\ge 1$, $T>0$, $D_1>0$, $D_2>0$ and $0<p\le q$ with $pq>1$.
Parabolic system~\eqref{eq:1.1} is one of the simplest parabolic systems
and it is an example of reaction-diffusion systems describing heat propagation
in a two component combustible mixture.
It has been studied extensively in many papers from various points of view,
see e.g. \cite{AHV, EH01, EH02, FLU, FI01, IKS, MST, QS, S, Z}.
See also \cite[Chapter~32]{QSBook} and references therein.
In this paper we prove the existence and the uniqueness
of the initial trace of a solution~$(u,v)$ to parabolic system \eqref{eq:1.1}
and study qualitative properties of the initial traces.
Furthermore, we obtain necessary conditions
on the existence of the solutions to problem~\eqref{eq:1.1} with
\begin{equation}
\label{eq:1.2}
(u(\cdot,0),v(\cdot,0))=(\mu,\nu)\quad\mbox{in}\quad{\bf R}^N,
\end{equation}
where $(\mu,\nu)$ is a pair of (nonnegative) Radon measures or nonnegative measurable functions in~${\bf R}^N$.

The study of the initial traces of the solutions to parabolic equations is a classical subject
and it has been investigated for various parabolic equations,
for example, the heat equation (see \cite{A, W}), the porous medium equation (see \cite{AC, BCP, HP}),
the parabolic $p$-Laplace equation (see \cite{DH, DH02}), the doubly nonlinear parabolic equation (see \cite{I, IKin, ZX}),
the fractional diffusion equation (see \cite{BSV}), the Finsler heat equation (see \cite{AIS})
and parabolic equations with nonlinear terms (see e.g.~\cite{ADi, BP, BCV, HI01, HI02, MV, Zao}).
For the semilinear heat equation
\begin{equation}
\label{eq:1.3}
\left\{
\begin{array}{ll}
\partial_t u=\Delta u+u^p & \quad\mbox{in}\quad{\bf R}^N\times(0,T),\vspace{3pt}\\
u\ge 0  & \quad\mbox{in}\quad{\bf R}^N\times(0,T),
\end{array}
\right.
\end{equation}
where $p>1$ and $T>0$,
by \cite{HI01} we have:
\begin{itemize}
 \item[{\rm (1)}]
 Let $u$ be a solution to \eqref{eq:1.3} for some $T>0$.
 Then there exists a unique (nonnegative) Radon measure $\mu$ in ${\bf R}^N$
 as the initial trace of $u$, that is,
 $$
 \underset{t\to+0}{\mbox{{\rm ess lim}}}
 \int_{{\bf R}^N}u(x,t)\varphi(x)\,dx=\int_{{\bf R}^N}\varphi(x)\,d\mu(x)
 $$
 for all $\varphi\in C_0({\bf R}^N)$.
 Furthermore, there exists a positive constant $\gamma_*$ depending only on $N$ and $p$
 such that
 \begin{equation}
 \label{eq:1.4}
 \sup_{x\in{\bf R}^N}\mu(B(x,\sigma))\le \gamma_*\,
 \sigma^{N-\frac{2}{p-1}}
 \end{equation}
 for $0<\sigma\le T^{\frac{1}{2}}$.
 In particular, in the case of  $p=1+2/N$,
  \begin{equation}
 \label{eq:1.5}
 \sup_{x\in{\bf R}^N}\mu(B(x,\sigma))\le \gamma_*\,
 \biggr[\log\biggr(e+\frac{T^{\frac{1}{2}}}{\sigma}\biggr)\biggr]^{-\frac{N}{2}}
 \end{equation}
 holds for $0<\sigma\le T^{\frac{1}{2}}$.
 \item[(2)]
 Let $u$ be a solution to \eqref{eq:1.3} for some $T>0$.
 Let $\mu$ be a (nonnegative) Radon measure in ${\bf R}^N$.
 Then $\mu$ is the initial trace of $u$ if and only if $u$ is a solution to \eqref{eq:1.3} with
 \begin{equation}
 \label{eq:1.6}
 u(\cdot,0)=\mu\quad\mbox{in}\quad{\bf R}^N.
 \end{equation}
\end{itemize}
See also \cite{ADi} and \cite{BP}.
Estimates~\eqref{eq:1.4} and \eqref{eq:1.5} give optimal necessary
conditions on the existence of solutions to problem~\eqref{eq:1.3} with \eqref{eq:1.6}
(see \cite[Theorems~1.3, 1.4 and 1.5]{HI01}).
Compared with semilinear parabolic equation~\eqref{eq:1.3},
much less is known about necessary conditions on the existence of solutions to semilinear parabolic systems
and there are no results corresponding to assertions~(1) and (2)
for parabolic system~\eqref{eq:1.1} even for the case of $D_1=D_2$.
\vspace{3pt}

In this paper we obtain necessary conditions on the existence of solutions to
semilinear parabolic system~\eqref{eq:1.1}
by being divided into the the following six cases (see Figure~1):
\begin{equation*}
\begin{split}
 & {\rm (A)}\quad\frac{q+1}{pq-1}<\frac{N}{2};\\
 & {\rm (B)}\quad \frac{q+1}{pq-1}=\frac{N}{2}\quad\mbox{and}\quad p<q ;
  \qquad\qquad\,\,
	{\rm (C)}\quad \frac{q+1}{pq-1}=\frac{N}{2}\quad\mbox{and}\quad p=q ;\\
 & {\rm (D)}\quad \frac{q+1}{pq-1}>\frac{N}{2}\quad\mbox{and}\quad q>1+\frac{2}{N};
	\qquad
	{\rm (E)}\quad \frac{q+1}{pq-1}>\frac{N}{2}\quad\mbox{and}\quad q=1+\frac{2}{N};\\
 & {\rm (F)}\quad \frac{q+1}{pq-1}>\frac{N}{2}\quad\mbox{and}\quad q<1+\frac{2}{N}.
\end{split}
\end{equation*}
\begin{figure}[htbp]
	\centering
	\subfloat{
		\begin{tikzpicture}[samples = 100]
			\draw[-{Latex}] (-0.5,0) -- (5.15,0) node[right] {\small $p$};
			\draw[-{Latex}] (0,-0.5) -- (0,5.15) node[above] {\small $q$};
			\path (0,0) node[above left] {\small O};
			\fill [ opacity = 0.2 ] plot [ domain = 0.2:1 ] ({\x}, {1/\x}) -- plot [ domain = 1:5 ] ({\x}, {\x});
			\draw plot[ domain = -0.4:5 ] ({\x}, {\x});
			\path (4,4) node[below right] {\small $p=q$};
      \draw [dashed] (3/5,5/3) -- (0,5/3) node[left] {\small $1+\frac{2}{N}$};
			\draw[color = white, thick] plot[ domain = 0.2:1 ] ({\x}, {1/\x});
			\draw[dashed] plot[ domain = 0.2:5 ] ({\x}, {1/\x});
			\draw[color = white, thick] plot[ domain = 1:5/3 ] ({\x}, {(5/3) / (\x - 2/3)});
			\draw[dashed] plot[ domain = 1:5/3 ] ({\x}, {(5/3) / (\x - 2/3)});
			\draw plot[ domain = 5/3:5 ] ({\x}, {(5/3) / (\x - 2/3)});
			\path (4.2,1.2) node {\small $\dfrac{q+1}{pq-1}=\dfrac{N}{2}$};
      \path (-0.8, 3.5) node {\small $pq=1$};
      \path (-0.8, 3.2) edge [-{Latex}, bend right] (1/2.7,2.7);
      \draw [thick, color = white] (3/5,5/3) -- (5/3,5/3);
			\draw [dashed] (3/5,5/3) -- (5/3,5/3);
      \draw [dashed] (5/3,5/3) -- (5/3,0) node[below] {\small $1+\frac{2}{N}$};
			\fill [color = white] (5/3, 5/3) circle (2pt);
			\draw (5/3, 5/3) circle (2pt);
			\fill [ color = white ] (1,1) circle (2pt);
			\draw (1,1) circle (2pt);
			\path (0.85,2.8) node {\small (D)};
			\path (2.3,3.3) node[above] {\small (A)};
			\path (1.03,1.1) node[above] {\small (F)};
		\end{tikzpicture}
	}
	\quad
	\subfloat{
		\begin{tikzpicture}[samples = 100]
			\draw[-{Latex}] (-0.5,0) -- (5.15,0) node[right] {\small $p$};
			\draw[-{Latex}] (0,-0.5) -- (0,5.15) node[above] {\small $q$};
			\path (0,0) node[above left] {\small O};
			\draw plot[ domain = -0.4:5 ] ({\x}, {\x});
			\path (4,4) node[below right] {\small $p=q$};
			\draw[dashed] plot[ domain = 0.2:5 ] ({\x}, {1/\x});
			\draw plot[ domain = 1:5 ] ({\x}, {(5/3) / (\x - 2/3)});
			\path (4.2,1.2) node {\small $\dfrac{q+1}{pq-1}=\dfrac{N}{2}$};
      \path (-0.8, 3.5) node {\small $pq=1$};
      \path (-0.8, 3.2) edge [-{Latex}, bend right] (1/2.7,2.7);
      \draw [dashed] (5/3,5/3) -- (5/3,0) node[below] {\small $1+\frac{2}{N}$};
      \draw [dashed] (3/5,5/3) -- (0,5/3) node[left] {\small $1+\frac{2}{N}$};
			\fill (5/3, 5/3) circle (2pt);
			\fill [ color = white ] (1,1) circle (2pt);
			\draw (1,1) circle (2pt);
			\draw (3/5,5/3) -- (5/3,5/3);
			\path (2.8,1.9) edge [-{Latex}, bend left = 30] (1.74,5/3);
			\path (2.84,1.8) node[above] {\small (C)};
			\path (0.8,2.3) edge [-{Latex}, bend right = 15] (1,1.7);
			\path (0.8,2.2) node[above] {\small (E)};
			\path (2,3.5) edge [-{Latex}, bend left] (1.25,3);
			\path (2,3.4) node[above] {\small (B)};
		\end{tikzpicture}
	}
	\caption{}\label{figure:1}
\end{figure}

\noindent
Our necessary conditions in cases~(A), (C) and  (F) can be regarded as generalizations of
the necessary conditions on the existence of solutions to problem~\eqref{eq:1.3} in the cases $p>1+2/N$,
$p=1+2/N$ and $1<p<1+2/N$, respectively.
On the other hand, our necessary conditions in cases (B), (D) and (E) are specific to parabolic systems.

We formulate the definition of the solution to \eqref{eq:1.1}.
For any Radon measure $\mu$ in $\mathbf{R}^N$,
we set
\begin{equation}
\label{eq:1.7}
[S(t)\mu](x):=\int_{\mathbf{R}^N} G(x-y,t)\, d\mu(y),
\quad\mbox{where}\quad
G(x,t)=(4\pi t)^{-\frac{N}{2}}\exp\left(-\frac{|x|^2}{4t}\right).
\end{equation}
We also write
$$
[S(t)\mu](x)=\int_{{\bf R}^N}G(x-y,t)\mu(y)\,dy
$$
if $\mu$ is a nonnegative measurable function in ${\bf R}^N$.
Then $S(t)\mu$ is the unique solution to the heat equation
$\partial_t u = \Delta u$ in ${\bf R}^N\times(0,\infty)$ with $u(0)=\mu$
under a suitable growth condition on the solution at the space infinity.
\begin{definition}
\label{Definition:1.1}
Let $0<T\le\infty$ and let $(u,v)$ be a pair of nonnegative measurable functions in ${\bf R}^N\times(0,T)$.
\vspace{3pt}
\newline
{\rm (i)} We say that $(u,v)$ is a solution to \eqref{eq:1.1} in ${\bf R}^N\times(0,T)$ if
$u(x,t)<\infty$, $v(x,t)<\infty$ and
\begin{equation*}
\begin{split}
  u(x,t) & =[S(D_1(t-\tau))u(\tau)](x)+\int_\tau^t [S_1(D_1(t-s))v(s)^p](x)\,ds,\\
  v(x,t) & =[S(D_2(t-\tau))v(\tau)](x)+\int_\tau^t [S_1(D_2(t-s))u(s)^q](x)\,ds,
\end{split}
\end{equation*}
for almost all $x\in{\bf R}^N$ and $0<\tau<t<T$.
\vspace{3pt}
\newline
{\rm (ii)} Let $\mu$ and $\nu$ be Radon measures in ${\bf R}^N$.
We say that $(u,v)$ is a solution to \eqref{eq:1.1} with \eqref{eq:1.2} in ${\bf R}^N\times(0,T)$ if
$u(x,t)<\infty$, $v(x,t)<\infty$ and
\begin{equation*}
\begin{split}
  u(x,t) & =[S(D_1(t))\mu](x)+\int_0^t [S_1(D_1(t-s))v(s)^p](x)\,ds,\\
  v(x,t) & =[S(D_2(t))\nu](x)+\int_0^t [S_1(D_2(t-s))u(s)^q](x)\,ds,
\end{split}
\end{equation*}
for almost all $x\in{\bf R}^N$ and $0<t<T$.
\end{definition}
\begin{remark}
\label{Remark:1.1}
Let $(u,v)$ be a solution to \eqref{eq:1.1} with \eqref{eq:1.2} in ${\bf R}^N\times(0,T)$
in the sense of Definition~{\rm \ref{Definition:1.1}}~{\rm (ii)},
where $T>0$. Then it follows from Fubini's theorem that
$(u,v)$ is a solution to \eqref{eq:1.1} in ${\bf R}^N\times(0,T)$ in the sense of Definition~{\rm \ref{Definition:1.1}}~{\rm (i)}.
\end{remark}

We are ready to state the main results of this paper.
Theorem~\ref{Theorem:1.1} shows the existence and the uniqueness of the initial trace
of a solution to \eqref{eq:1.1}
and gives an upper estimate of the strength
of the singularity of the initial trace.
\begin{theorem}
\label{Theorem:1.1}
Let $N\ge 1$, $D_1>0$, $D_2>0$, $0<p\le q$ with $pq>1$ and $T>0$.
Let $(u,v)$ be a solution to \eqref{eq:1.1} in ${\bf R}^N\times(0,T)$.
Then there exists a unique pair of Radon measures $(\mu,\nu)$ such that
\begin{equation}
\label{eq:1.8}
\begin{split}
 & \underset{t\to+0}{\mbox{{\rm ess lim}}}\,
\int_{{\bf R}^N}u(y,t)\phi(y)\,dy=\int_{{\bf R}^N}\phi(y)\,d\mu(y),\\
 & \underset{t\to+0}{\mbox{{\rm ess lim}}}\,
\int_{{\bf R}^N}v(y,t)\phi(y)\,dy=\int_{{\bf R}^N}\phi(y)\,d\nu(y),
\end{split}
\end{equation}
for all $\phi\in C_0({\bf R}^N)$.
Furthermore, there exists a positive constant $\gamma$ depending only on $N$, $D_1$, $D_2$, $p$ and $q$ such that
\begin{equation}
\label{eq:1.9}
\sup_{x\in{\bf R}^N}\,\mu(B(x,\sigma))\le\gamma \sigma^{N-\frac{2(p+1)}{pq-1}},
\quad
\sup_{x\in{\bf R}^N}\,\nu(B(x,\sigma))\le\gamma \sigma^{N-\frac{2(q+1)}{pq-1}},
\end{equation}
for $0<\sigma\le T^\frac{1}{2}$.
\end{theorem}
\begin{remark}
\label{Remark:1.2}
In case~{\rm (F)}, it follows from $p\le q$ and $q<1+2/N$ that
$$
\frac{p+1}{pq-1}-\frac{N}{2}>\frac{p+1}{pq-1}-\frac{1}{q-1}
=\frac{(p+1)(q-1)-(pq-1)}{(pq-1)(q-1)}
=\frac{q-p}{(pq-1)(q-1)}\ge 0.
$$
Then we see that the functions
$$
(0,\infty)\ni\sigma\mapsto \sigma^{N-\frac{2(p+1)}{pq-1}}\quad\mbox{and}\quad
(0,\infty)\ni\sigma\mapsto \sigma^{N-\frac{2(q+1)}{pq-1}}
$$
are monotonically decreasing. This implies that
estimate~\eqref{eq:1.9} is equivalent to
$$
\sup_{x\in{\bf R}^N}\,\mu(B(x,T^{\frac{1}{2}}))\le\gamma T^{\frac{N}{2}-\frac{p+1}{pq-1}},
\quad
\sup_{x\in{\bf R}^N}\,\nu(B(x,T^{\frac{1}{2}}))\le\gamma T^{\frac{N}{2}-\frac{q+1}{pq-1}}.
$$
Similarly, in cases {\rm (D)} and {\rm (E)},
the estimates for $\nu$ of \eqref{eq:1.9} are equivalent to
\begin{equation*}
  \sup_{x\in{\bf R}^N}\,\nu(B(x,T^{\frac{1}{2}}))\le\gamma T^{\frac{N}{2}-\frac{q+1}{pq-1}}.
\end{equation*}
\end{remark}
As an application of Theorem~\ref{Theorem:1.1},
we obtain the following result, which corresponds to assertion~(2).
\begin{theorem}
\label{Theorem:1.2}
Assume the same conditions as in Theorem~{\rm\ref{Theorem:1.1}}.
Let $\mu$ and $\nu$ be Radon measures in ${\bf R}^N$ and $0<T\le\infty$.
\begin{itemize}
  \item[{\rm (i)}]
  Let $(u,v)$ be a solution to \eqref{eq:1.1} with \eqref{eq:1.2} in ${\bf R}^N\times(0,T)$.
  Then \eqref{eq:1.8} holds.
  \item[{\rm (ii)}]
  Let $(u,v)$ be a solution to \eqref{eq:1.1} in ${\bf R}^N\times (0,T)$.
  Assume \eqref{eq:1.8}.
  Then $(u,v)$ is a solution to \eqref{eq:1.1} with \eqref{eq:1.2} in ${\bf R}^N\times(0,T)$.
\end{itemize}
\end{theorem}
In Theorem~\ref{Theorem:1.3} we improve estimate~\eqref{eq:1.9} in cases~(B), (C), (D) and (E).
\begin{theorem}
\label{Theorem:1.3}
Assume the same conditions as in Theorem~{\rm\ref{Theorem:1.1}}.
Let $(\mu,\nu)$ be a pair of Radon measures in ${\bf R}^N$ satisfying \eqref{eq:1.8}.
\begin{itemize}
  \item[{\rm (1)}]
  Consider case~{\rm (B)}.
  Then there exists $\gamma_B>0$ such that
  $$
  \sup_{x\in{\bf R}^N}\,\int_0^\sigma \tau^{-1}\left[\frac{\mu(B(x,\tau))}{\tau^{N-\frac{2(p+1)}{pq-1}}}\right]^q\,d\tau
  +\sup_{x\in{\bf R}^N}\,\nu(B(x,\sigma))
  \le\gamma_B\biggr[\log\left(e+\frac{T^\frac{1}{2}}{\sigma}\right)\biggr]^{-\frac{1}{pq-1}}
  $$
  for $0<\sigma\le T^\frac{1}{2}$.
  \item[{\rm (2)}]
  Consider case~{\rm (C)}.
  Then there exists $\gamma_C>0$ such that
  $$
  \sup_{x\in{\bf R}^N}\,\mu(B(x,\sigma))+\sup_{x\in{\bf R}^N}\,\nu(B(x,\sigma))
  \le\gamma_C \left[ \log \left( e + \frac{T^\frac{1}{2}}{\sigma} \right) \right]^{-\frac{N}{2}}
  $$
  for $0<\sigma\le T^\frac{1}{2}$.
  \item[{\rm (3)}]
  Consider case~{\rm (D)}.
  Then there exists $\gamma_D>0$ such that
  $$
  \sup_{x\in{\bf R}^N}\,\int_0^{T^\frac{1}{2}}
  \tau^{-1}\left[\frac{\mu(B(x,\tau))}{\tau^{N-\frac{N+2}{q}}}\right]^q\,d\tau
  +\sup_{x\in{\bf R}^N}\,\nu(B(x,T^\frac{1}{2})) \le \gamma_D T^{\frac{N}{2} - \frac{q+1}{pq-1}}.
  $$
  \item[{\rm (4)}]
  Consider case~{\rm (E)}.
  Then there exists $\gamma_E>0$ such that
  $$
  \sup_{x\in{\bf R}^N}\,\int_0^{T^\frac{1}{2}} \tau^{-1}\mu(B(x,\tau))^q\,d\tau
  +\sup_{x\in{\bf R}^N}\,\nu(B(x,T^\frac{1}{2})) \le\gamma_E T^{\frac{N}{2} - \frac{q+1}{pq-1}}.
  $$
\end{itemize}
Here the constants $\gamma_B$, $\gamma_C$, $\gamma_D$ and $\gamma_E$ depend only on
$N$, $D_1$, $D_2$, $p$ and $q$.
\end{theorem}
\begin{remark}
  \label{Remark:1.3}
  As in Remark~{\rm \ref{Remark:1.2}},
  the estimates in assertions {\rm (3)} and {\rm (4)} are equivalent to
  \begin{gather*}
    \sup_{x\in{\bf R}^N}\,\int_0^{\sigma}
    \tau^{-1}\left[\frac{\mu(B(x,\tau))}{\tau^{N-\frac{N+2}{q}}}\right]^q\,d\tau
    +\sup_{x\in{\bf R}^N}\,\nu(B(x,\sigma)) \le \gamma_D \sigma^{N - \frac{2(q+1)}{pq-1}},
    \\
    \sup_{x\in{\bf R}^N}\,\int_0^{\sigma} \tau^{-1}\mu(B(x,\tau))^q\,d\tau
    +\sup_{x\in{\bf R}^N}\,\nu(B(x,\sigma)) \le\gamma_E \sigma^{N - \frac{2(q+1)}{pq-1}},
  \end{gather*}
  for all $0<\sigma\le T^{1/2}$, respectively.
\end{remark}
As corollaries of our theorems,
we obtain the following results, which show the validity
of estimate~\eqref{eq:1.9} and Theorem~\ref{Theorem:1.3}.
Note that Corollary~\ref{Corollary:1.1} has been already proved
in \cite{EH02} for $D_1=D_2$ and in \cite{FLU} for $D_1\not=D_2$.
\begin{corollary}
\label{Corollary:1.1}
Let $N\ge 1$, $0<p\le q$ with $pq>1$. Assume that
\begin{equation}
\label{eq:1.10}
\frac{q+1}{pq-1}\ge\frac{N}{2}.
\end{equation}
Then problem~\eqref{eq:1.1} possesses no global-in-time nontrivial solutions.
\end{corollary}
\begin{remark}
\label{Remark:1.4}
Let $N\ge 1$ and $0<p\le q$ with $pq>1$.
Then \eqref{eq:1.10} is an optimal condition for the nonexistence of global-in-time
nontrivial solutions to \eqref{eq:1.1}.
Indeed, if
$$
\frac{q+1}{pq-1}<\frac{N}{2},
$$
then problem~\eqref{eq:1.1} possesses a global-in-time positive solution.
See {\rm\cite{EH02}} and {\rm\cite{FLU}}.
\end{remark}
\begin{corollary}
\label{Corollary:1.2}
Let $N\ge 1$ and $0<p\le q$ with $pq>1$.
\begin{itemize}
\item[{\rm (a)}]
Consider case {\rm (A)}.
Let $(\mu,\nu)$ be a pair of nonnegative measurable functions in ${\bf R}^N$ such that
\begin{equation}
\label{eq:1.11}
\mu(x)\ge c_{a,1} |x|^{-\frac{2(p+1)}{pq-1}},\qquad
\nu(x)\ge c_{a,2} |x|^{-\frac{2(q+1)}{pq-1}},
\end{equation}
in a neighborhood of $0$,
where $c_{a,1}$, $c_{a,2}>0$.
Then
problem~\eqref{eq:1.1} with \eqref{eq:1.2} possesses no local-in-time solutions
if either $c_{a,1}$ or $c_{a,2}$ is sufficiently large.
\item[{\rm (b)}]
Consider case {\rm (B)}.
Let $(\mu,\nu)$ be a pair of nonnegative measurable functions in ${\bf R}^N$ such that
\begin{equation}
\label{eq:1.12}
\mu(x)\ge c_{b,1} |x|^{-\frac{2(p+1)}{pq-1}}|\log|x||^{-\frac{p}{pq-1}},
\qquad
\nu(x)\ge c_{b,2} |x|^{-N}|\log|x||^{-\frac{1}{pq-1}-1},
\end{equation}
in a neighborhood of $0$, where $c_{b,1}$, $c_{b,2}>0$.
Then problem~\eqref{eq:1.1} with \eqref{eq:1.2} possesses no local-in-time solutions
if either $c_{b,1}$ or $c_{b,2}$ is sufficiently large.
\item[{\rm (c)}]
Consider case {\rm (C)}.
Let $(\mu,\nu)$ be a pair of nonnegative measurable functions in ${\bf R}^N$ such that
\begin{equation}
\label{eq:1.13}
\mu(x)\ge c_{c,1}|x|^{-N}|\log|x||^{-\frac{N}{2}-1},
\qquad
\nu(x)\ge c_{c,2}|x|^{-N}|\log|x||^{-\frac{N}{2}-1},
\end{equation}
in a neighborhood of $0$, where $c_{c,1}$, $c_{c,2}>0$.
Then problem~\eqref{eq:1.1} with \eqref{eq:1.2} possesses no local-in-time solutions
if either $c_{c,1}$ or $c_{c,2}$ is sufficiently large.
\item[{\rm (d)}]
Consider case {\rm (D)}.
Let $\mu$ be a nonnegative measurable function in ${\bf R}^N$ such that
\begin{equation}
\label{eq:1.14}
\mu(x)\ge |x|^{-\frac{N+2}{q}}h_1(|x|)
\end{equation}
in a neighborhood of $0$.
Here $h_1$ is a positive continuous function in $(0,1]$ such that
$s^{-\epsilon}h(s)$ is monotonically decreasing in $(0,\delta)$ for some $\epsilon>0$ and $\delta>0$.
Let $\nu$ be a Radon measure in ${\bf R}^N$.
Assume either
\begin{equation}
\label{eq:1.15}
\int_0^1 \tau^{-1}h_1(\tau)^q\,d\tau=\infty\quad\mbox{or}\quad
\sup_{x\in{\bf R}^N}\,\nu(B(x,1))=\infty.
\end{equation}
Then problem~\eqref{eq:1.1} with \eqref{eq:1.2} possesses no local-in-time solutions.
\item[{\rm (e)}]
Consider case {\rm (E)}.
Let $\mu$ be a nonnegative measurable function in ${\bf R}^N$
such that
\begin{equation}
\label{eq:1.16}
\mu(x)\ge |x|^{-N}h_2(|x|)
\end{equation}
in a neighborhood of $0$, where $h_2$ is a positive continuous function in $(0,1]$.
Let $\nu$ be a Radon measure in ${\bf R}^N$.
Assume either
$$
\int_0^1\left[\int_0^r \tau^{-1}h_2(\tau)\,d\tau\right]^qr^{-1}\,dr=\infty
\quad\mbox{or}\quad
\sup_{x\in{\bf R}^N}\,\nu(B(x,1))=\infty.
$$
Then problem~\eqref{eq:1.1} with \eqref{eq:1.2} possesses no local-in-time solutions.
\item[{\rm (f)}]
Consider case {\rm (F)}.
Let $(\mu,\nu)$ be a pair of Radon measures in ${\bf R}^N$.
Then problem~\eqref{eq:1.1} with \eqref{eq:1.2} possesses no local-in-time solutions
if either
$$
\sup_{x\in{\bf R}^N}\mu(B(x,1))=\infty\quad\mbox{or}\quad \sup_{x\in{\bf R}^N}\nu(B(x,1))=\infty.
$$
\end{itemize}
\end{corollary}
Corollary~\ref{Corollary:1.2} actually gives optimal conditions
for the nonexistence of local-in-time solutions to problem~\eqref{eq:1.1} with \eqref{eq:1.2}.
Indeed, in case~(A),
if
$$
0\le\mu(x)\le c_{a,1}'|x|^{-\frac{2(p+1)}{pq-1}}
\quad\mbox{and}\quad
0\le\nu(x)\le c_{a,1}'|x|^{-\frac{2(q+1)}{pq-1}}
\quad\mbox{in}\quad{\bf R}^N
$$
hold with sufficiently small positive constants $c_{a,1}'$ and $c_{a,2}'$,
then problem~\eqref{eq:1.1} with \eqref{eq:1.2} possesses a global-in-time solution
(see \cite[Corollaries~3.1 and 3.2]{IKS}).
Similar results also hold in the other cases~(B), (C), (D), (E) and (F).
See a forthcoming paper~\cite{FI02}.
\vspace{3pt}

The rest of this paper is organized as follows.
In Section~2 we recall some properties of the Gauss kernel and prove some preliminary lemmas.
In Section~3 we modify the arguments in \cite{EH02} to prove Theorem~\ref{Theorem:1.1}.
In Section~4 we follow the arguments in \cite{HI01} to prove Theorem~\ref{Theorem:1.2}.
In Sections~5, 6 and 7 we develop the arguments in \cite{HI01} to prove Theorem~\ref{Theorem:1.3}.
The proof of Theorem~\ref{Theorem:1.3} requires more delicate arguments than those in \cite{HI01}.
Section~8 is devoted to the proofs of Corollaries~\ref{Corollary:1.1} and \ref{Corollary:1.2}.
\section{Preliminary lemmas}
In this section we recall some properties of $S(t)\mu$ and prove some preliminary lemmas.
Throughout this paper, by the letter $C$
we denote generic positive constants
depending only on $N$, $D_1$, $D_2$, $p$ and $q$,
and they may have different values also within the same line.

It follows from \eqref{eq:1.7} that
$$
\int_{{\bf R}^N}G(x,t)\,dx=1.
$$
Then Jensen's inequality implies that
\begin{equation}
\label{eq:2.1}
\begin{split}
[S(t)\mu](x) & =\int_{{\bf R}^N}G(x-y,t)\mu(y)\,dy\\
 & \le\left(\int_{{\bf R}^N}G(x-y,t)\mu(y)^\alpha\,dy\right)^{\frac{1}{\alpha}}
=[S(t)\mu^\alpha](x)^{\frac{1}{\alpha}},
\quad x\in{\bf R}^N,\,\,t>0,
\end{split}
\end{equation}
for any $\alpha\ge 1$.
On the other hand, since $G=G(x,t)$ satisfies
\begin{equation}
\label{eq:2.2}
G(x,t)=\int_{{\bf R}^N}G(x-y,t-s)G(y,s)dy,\quad x\in{\bf R}^N,\,\,0<s<t,
\end{equation}
we find
$$
[S(t)\mu](x)=[S(t-s)(S(s)\mu)](x),\quad x\in{\bf R}^N,\,\,0<s<t.
$$
Furthermore, we have the following lemmas.
\begin{lemma}
\label{Lemma:2.1}
{\rm (i)} Let $\mu\in L^1_{{\rm loc}}({\bf R}^N)$ be nonnegative in ${\bf R}^N$.
Then there exists $C>0$ depending only on $N$ such that
$$
\int_{B(0,\rho)}G(x-y,t)\mu(y)\,dy\ge CG\left(x,\frac{t}{2}\right)\int_{B(0,\rho)}\mu(y)\,dy
$$
for $x\in{\bf R}^N$, $\rho>0$ and $t\ge\rho^2$.
\vspace{3pt}
\newline
{\rm (ii)} Let $\mu$ be a Radon measure in ${\bf R}^N$.
Then
$$
\int_{B(0,\rho)}G(x-y,t)\,d\mu(y)\ge CG\left(x,\frac{t}{2}\right)\int_{B(0,\rho)}\,d\mu(y)
$$
for $x\in{\bf R}^N$, $\rho>0$ and $t\ge\rho^2$.
\end{lemma}
\textbf{Proof.}
We prove assertion~(i).
Let $x\in{\bf R}^N$, $\rho>0$ and $t\ge\rho^2$.
Since $|x-y|^2 \le 2(|x|^2+|y|^2)$,
for $y\in B(0,\rho)$,
we have
\begin{equation*}
  \exp\left( -\frac{|x-y|^2}{4t} \right)
  \ge \exp\left( -\frac{|x|^2}{2t} \right)
  \cdot \exp\left( -\frac{|y|^2}{2t} \right)
  \ge e^{-\frac{1}{2}}\exp\left( -\frac{|x|^2}{2t} \right).
\end{equation*}
This implies that
\begin{equation*}
  \int_{B(0,\rho)} G(x-y,t)\mu(y)\, dy
  \ge 2^{-\frac{N}{2}}e^{-\frac{1}{2}} G\left(x,\frac{t}{2}\right)\int_{B(0,\rho)}\mu(y)\,dy,
\end{equation*}
that is, assertion~(i) follows. Similarly, assertion~(ii) follows.
The proof is complete.
$\Box$\vspace{5pt}
\begin{lemma}
\label{Lemma:2.2}
Let $\alpha>0$, $\beta>0$ and $L\ge 0$. Then
\begin{equation*}
  \begin{split}
    &
    \int_{{\bf R}^N}G(x-y,t-s)G\left(y,\frac{s+L}{\alpha}\right)^\beta\,dy\\
    &  \ge\left(\frac{4\pi}{\alpha}\right)^{-\frac{N}{2}(\beta-1)}\beta^{-\frac{N}{2}}
    \left(\frac{\min\{\alpha\beta,1\}}{\max\{\alpha\beta,1\}}\right)^{\frac{N}{2}}
    (s+L)^{-\frac{N}{2}(\beta-1)}G\left(x,\frac{t+L}{\max\{\alpha\beta,1\}}\right)
  \end{split}
\end{equation*}
for $x\in{\bf R}^N$ and $0<s<t$.
\end{lemma}
{\bf Proof.}
Let $x$, $y\in{\bf R}^N$ and $0<s<t$.
It follows that
\begin{equation*}
  \begin{aligned}
G\left(y,\frac{s+L}{\alpha}\right)^\beta
    & =
    \left[\frac{4\pi(s+L)}{\alpha}\right]^{-\frac{N(\beta-1)}{2}}
    \left[\frac{4\pi(s+L)}{\alpha\beta}\right]^{-\frac{N}{2}}\beta^{-\frac{N}{2}}
    \exp\left(-\frac{\alpha \beta |y|^2}{4(s+L)}\right)
    \\
    & =
    \left( \frac{4\pi}{\alpha} \right)^{-\frac{N(\beta-1)}{2}}
    \beta^{-\frac{N}{2}}(s+L)^{-\frac{N}{2}(\beta-1)}
    G\left(y,\frac{s+L}{\alpha\beta}\right).
  \end{aligned}
\end{equation*}
This together with \eqref{eq:2.2} implies that
\begin{equation}
\label{eq:2.3}
\begin{split}
  & \int_{{\bf R}^N}G(x-y,t-s)G\left(y,\frac{s+L}{\alpha}\right)^\beta\,dy \\
  & =\left( \frac{4\pi}{\alpha} \right)^{-\frac{N(\beta-1)}{2}}
    \beta^{-\frac{N}{2}}(s+L)^{-\frac{N}{2}(\beta-1)}
    \int_{{\bf R}^N}G(x-y,t-s)G\left(y,\frac{s+L}{\alpha\beta}\right)\,dy\\
  & =\left(\frac{4\pi}{\alpha}\right)^{-\frac{N}{2}(\beta-1)}\beta^{-\frac{N}{2}}
     (s+L)^{-\frac{N}{2}(\beta-1)}G\left(x,t-s+\frac{s+L}{\alpha\beta}\right).
\end{split}
\end{equation}
Since
$$
\frac{t+L}{\max\{1,\alpha\beta\}}
\le t-s+\frac{s+L}{\alpha\beta}
\le \frac{t+L}{\min\{1,\alpha\beta\}},
$$
by \eqref{eq:2.3} we obtain the desired inequality,
and the proof is complete.
$\Box$
\vspace{5pt}

We prove the following basic lemma.
\begin{lemma}
\label{Lemma:2.3}
For any $k\ge 1$, there exists $m\in\{1,2,\dots\}$ with the following property:
For any $R>0$ and $z\in{\bf R}^N$, there exists $\{z_i\}_{i=1}^m\subset{\bf R}^N$ such that
\begin{equation}
\label{eq:2.4}
B(z,kR)\subset\bigcup_{i=1}^m B(z_i,R).
\end{equation}
\end{lemma}
{\bf Proof.}
For any $k\ge 1$,
we can find $m\in\{1,2,\dots\}$ and $\{x_i\}_{i=1}^m\subset{\bf R}^N$ such that
$$
B(0,k)\subset\bigcup_{i=1}^m B(x_i,1).
$$
Then \eqref{eq:2.4} holds with $z_i:=z+Rx_i$ $(i=1,\dots,m)$, and
Lemma~\ref{Lemma:2.3} follows.
$\Box$
\vspace{3pt}\newline
As a direct consequence of Lemma~\ref{Lemma:2.3},
we have:
\begin{lemma}
\label{Lemma:2.4}
Let $\mu$ be a Radon measure in ${\bf R}^N$ such that
$$
\sup_{x\in{\bf R}^N}\mu(B(x,R))<\infty
$$
for some $R>0$. Then, for any $k\ge 1$, there exists $m\in\{1,2,\dots\}$ such that
$$
\sup_{x\in{\bf R}^N}\mu(B(x,kR))\le m\sup_{x\in{\bf R}^N}\mu(B(x,R))<\infty.
$$
\end{lemma}

We prepare the Besicovitch covering theorem, which will be used in the proof of Theorem~\ref{Theorem:1.2}.
See \cite{EG}.
\begin{lemma}
  There exists an integer $m$ depending only on $N$ with the following property:
  If $\mathcal{F}$ is a collection of closed balls in $\mathbf{R}^N$ with positive radius such that
  \begin{equation*}
    \sup\left\{ \operatorname{diam} B :\, B \in \mathcal{F} \right\} < \infty
  \end{equation*}
  and $A$ is the set of centers of balls in $\mathcal{F}$,
  then there exist $\mathcal{B}_1,\ldots,\mathcal{B}_{m}\subset \mathcal{F}$ such that
  $\mathcal{B}_{i}$ $(i=1,\ldots,m)$ is a countable collection of disjoint balls in $\mathcal{F}$ and
  \begin{equation*}
    A \subset \bigcup_{i=1}^m \bigcup_{B\in \mathcal{B}_i} B.
  \end{equation*}
\end{lemma}
\section{Proof of Theorem~\ref{Theorem:1.1}}
We prove Theorem~\ref{Theorem:1.1}.
For this aim, it suffices to consider the case of $T=1$ and $D:=\min\{D_1,D_2\}=1$.
Indeed, let $(u,v)$ be a solution to \eqref{eq:1.1} in $\mathbf{R}^N \times(0,T)$, where $T>0$.
Let $\lambda>0$ be such that $\lambda^2\min\{D_1,D_2\}=1$. Set
$$
\tilde{u}(x,t) := T^{\frac{p+1}{pq-1}} u(\lambda^{-1}T^{\frac{1}{2}}x, Tt),
\qquad
\tilde{v}(x,t) := T^{\frac{q+1}{pq-1}} u(\lambda^{-1}T^{\frac{1}{2}}x, Tt),
$$
in $\mathbf{R}^N \times (0,1)$.
Then $(\tilde{u}, \tilde{v})$ satisfies
\begin{equation*}
\left\{
\begin{array}{ll}
\partial_t \tilde{u}=\tilde{D}_1\Delta \tilde{u}+\tilde{v}^p & \quad\mbox{in}\quad{\bf R}^N\times(0,1),\vspace{3pt}\\
\partial_t \tilde{v}=\tilde{D}_2\Delta \tilde{v}+\tilde{u}^q & \quad\mbox{in}\quad{\bf R}^N\times(0,1),\vspace{3pt}\\
\tilde{u},\tilde{v}\ge 0 & \quad\mbox{in}\quad{\bf R}^N\times(0,1),
\end{array}
\right.
\end{equation*}
where $\tilde{D}_1=\lambda^2 D_1$ and $\tilde{D}_2=\lambda^2 D_2$.
Here $\min\{\tilde{D}_1,\tilde{D}_2\}=1$.
\vspace{3pt}

Let $(u,v)$ be a solution to \eqref{eq:1.1} in $\mathbf{R}^N \times(0,1)$ and $D=1$.
Since
\begin{equation}
\label{eq:3.1}
G(x,D_it)=(4\pi D_i)^{-\frac{N}{2}}\exp\left(-\frac{|x|^2}{4D_i t}\right)
\ge(D')^{-\frac{N}{2}}G(x,t),\qquad x\in{\bf R}^N,\,\,t>0,
\end{equation}
where $i=1,2$ and $D':=\max\{D_1,D_2\}$,
by Definition~\ref{Definition:1.1}~(i) we see that
\begin{equation}
\label{eq:3.2}
\begin{split}
\infty>(D')^{\frac{N}{2}}u(x,t) & \ge [S(t-\tau)u(\tau)](x)+\int_\tau^t [S(t-s)v(s)^p](x)\,ds,\\
\infty>(D')^{\frac{N}{2}}v(x,t) & \ge [S(t-\tau)v(\tau)](x)+\int_\tau^t [S(t-s)u(s)^q](x)\,ds,
\end{split}
\end{equation}
for almost all $x\in{\bf R}^N$ and $0<\tau<t<1$.
Then we apply the same argument as in the proof of \cite[Lemma~2.3]{HI01}
and we can prove the existence and the uniqueness of the initial trace of $(u,v)$.
So it remains to prove \eqref{eq:1.9} for the proof of Theorem~\ref{Theorem:1.1}.
\vspace{5pt}

Let
\begin{equation}
\label{eq:3.3}
z\in{\bf R}^N,\qquad
0<\rho<1/\sqrt{5},\qquad
0<\tau\le\rho^2<1.
\end{equation}
Then the functions
\begin{equation}
\label{eq:3.4}
U(x,t):=u(x+z,t+\tau+2\rho^2),
\qquad
V(x,t):=v(x+z,t+\tau+2\rho^2),
\end{equation}
are defined for almost all $x\in{\bf R}^N$ and $t\in[0,t_*)$, where
$t_*:=1-\tau-2\rho^2\ge 2\rho^2$.
By \eqref{eq:3.2} we have
\begin{equation}
\label{eq:3.5}
\begin{split}
\infty>(D')^{\frac{N}{2}}U(x,t) & \ge [S(t)U(0)](x)+\int_0^t [S(t-s)V(s)^p](x)\,ds,\\
\infty>(D')^{\frac{N}{2}}V(x,t) & \ge [S(t)V(0)](x)+\int_0^t [S(t-s)U(s)^q](x)\,ds,
\end{split}
\end{equation}
for almost all $x\in{\bf R}^N$, $t\in[0,t_*)$ and $\tau\in(0,\rho^2]$.
By Lemma~\ref{Lemma:2.1} and \eqref{eq:3.2}
we see that
\begin{equation}
\label{eq:3.6}
\begin{split}
U(x,0) & =u(x+z,\tau+2\rho^2)
\ge (D')^{-\frac{N}{2}}[S(2\rho^2)u(\tau)](x+z)\\
 & \ge (D')^{-\frac{N}{2}}\int_{B(z,\rho)}G(x+z-y,2\rho^2)u(y,\tau)\,dy
\ge CM_{z,\tau}(u)G(x,\rho^2),\\
V(x,0) & =v(x+z,\tau+2\rho^2)
\ge (D')^{-\frac{N}{2}}[S(2\rho^2)v(\tau)](x+z)\\
 & \ge (D')^{-\frac{N}{2}}\int_{B(z,\rho)}G(x+z-y,2\rho^2)v(y,\tau)\,dy
\ge CM_{z,\tau}(v)G(x,\rho^2),
\end{split}
\end{equation}
for almost all $x\in \mathbf{R}^N$ and $\tau\in(0,\rho^2]$, where
$$
M_{z,\tau}(u):=\int_{B(z,\rho)} u(x,\tau)\,dx,
\qquad
M_{z,\tau}(v):=\int_{B(z,\rho)} v(x,\tau)\,dx.
$$
By \eqref{eq:3.5} and \eqref{eq:3.6} we have
\begin{equation}
\label{eq:3.7}
U(x,t)\ge c_*M_{z,\tau}(u)G(x,t+\rho^2),\quad
V(x,t)\ge c_*M_{z,\tau}(v)G(x,t+\rho^2),
\end{equation}
for almost all $x\in{\bf R}^N$, $t\in[0,t_*)$ and $\tau\in(0,\rho^2]$.
Here $c_*$ is a positive constant depending only on $N$, $p$, $q$, $D_1$ and $D_2$.
\begin{lemma}
\label{Lemma:3.1}
Let $\rho\in(0,1)$ and $\tau\in(0,1)$ and assume \eqref{eq:3.3}.
Furthermore, assume that there exist positive constants $a$, $b$ and $c\ge 1$ such that
\begin{equation}
\label{eq:3.8}
U(x,t)\ge at^bG(x,t+\rho^2)^c
\end{equation}
for almost all $x\in{\bf R}^N$ and $t\in[0,2\rho^2]$.
Then there exists a positive constant $\gamma_1$ depending only on $N$, $D_1$, $D_2$, $p$ and $q$
such that
$$
U(x,t)\ge \gamma_1\frac{a^{pq}}{(qb+1)^p(pqb+p+1)}t^{pqb+p+1}G(x,t+\rho^2)^{pqc}
$$
for almost all $x\in{\bf R}^N$ and $t\in[0,2\rho^2]$.
\end{lemma}
{\bf Proof.}
We follow the argument in the proof of \cite[Lemma~2.4]{EH02} to prove Lemma~\ref{Lemma:3.1}.
Since $q>1$,
by \eqref{eq:2.1}, \eqref{eq:2.2}, \eqref{eq:3.5} and \eqref{eq:3.8}
we have
\begin{equation}
  \label{eq:3.9}
  \begin{split}
    V(x,t)  & \ge
    (D')^{-\frac{N}{2}}\int_0^t [S(t-s)U(s)^q](x)\,ds \\
    & \ge Ca^q
    \int_0^t \int_{\mathbf{R}^N} s^{qb}G(x-y,t-s)G(y,s+\rho^2)^{qc} \, dy\, ds\\
    & \ge Ca^q
    \int_0^t s^{qb}\left(\int_{\mathbf{R}^N} G(x-y,t-s)G(y,s+\rho^2)\, dy\right)^{qc} \, ds\\
    & \ge C\frac{a^q}{qb+1} t^{qb+1}G(x,t+\rho^2)^{qc}
  \end{split}
\end{equation}
for almost all $x\in{\bf R}^N$ and $t\in[0,2\rho^2]$.
Since $pq>1$,
repeating the above argument
and using \eqref{eq:3.9}, instead of \eqref{eq:3.8},
we obtain
\begin{equation*}
  \begin{split}
    U(x,t)  & \ge
    (D')^{-\frac{N}{2}}\int_0^t [S(t-s)V(s)^p](x)\,ds \\
    & \ge C\frac{a^{pq}}{(qb+1)^p}
    \int_0^t \int_{\mathbf{R}^N} s^{p(qb+1)}G(x-y,t-s)G(y,s+\rho^2)^{pqc} \, dy\, ds\\
    & \ge C
    \frac{a^{pq}}{(qb+1)^p}
    \int_0^t s^{p(qb+1)}
    \left(\int_{\mathbf{R}^N} G(x-y,t-s)G(y,s+\rho^2)\, dy\right)^{pqc} \, ds\\
    & \ge C\frac{a^{pq}}{(qb+1)^p(pqb+p+1)}t^{pqb+p+1}G(x,t+\rho^2)^{pqc}
  \end{split}
\end{equation*}
for almost all $x\in{\bf R}^N$ and $t\in[0,2\rho^2]$.
Thus Lemma~\ref{Lemma:3.1} follows.
$\Box$
\vspace{5pt}\newline
Combining Lemma~\ref{Lemma:3.1} and \eqref{eq:3.7},
we have:
\begin{lemma}
\label{Lemma:3.2}
Let $(u,v)$ satisfy \eqref{eq:3.2}.
Then there exists $C_*>0$ such that
$$
\sup_{z\in{\bf R}^N}\int_{B(z,\rho)}u(y,\tau)\,dy\le C_*\rho^{N-\frac{2(p+1)}{pq-1}}
$$
for all $\rho\in(0,1/\sqrt{5})$ and almost all $\tau\in(0,\rho^2]$.
\end{lemma}
{\bf Proof.}
Let $z\in{\bf R}^N$ and $\rho\in(0,1/\sqrt{5})$.
Combining Lemma~\ref{Lemma:3.1} and \eqref{eq:3.7}, we obtain
\begin{equation}
\label{eq:3.10}
\infty>U(x,t)\ge a_n t^{b_n} G(x,t+\rho^2)^{c_n},
\quad n=0,1,2,\dots,
\end{equation}
for almost all $x\in{\bf R}^N$, $t\in[0,2\rho^2]$ and $\tau\in(0,\rho^2]$.
Here $\{a_n\}_{n=0}^\infty$, $\{b_n\}_{n=0}^\infty$ and $\{c_n\}_{n=0}^\infty$ are sequences defined by
$$
a_{n+1} :=\gamma_1\frac{a_n^{pq}}{(qb_n+1)^p(pqb_n+p+1)},\quad
b_{n+1} := pqb_n+p+1,\quad
c_{n+1} := pq c_n,
$$
for $n=0,1,2,\dots$ with
$$
a_0 := c_*M_{z,\tau}(u),\quad b_0 :=0, \quad c_0 :=1.
$$
Here $c_*$ and $\gamma_1$ are as in \eqref{eq:3.7} and Lemma~\ref{Lemma:3.1}, respectively.
Then we observe that
\begin{equation}
\label{eq:3.11}
b_n =[(pq)^n-1]\frac{p+1}{pq-1}, \qquad
c_n =(pq)^n,
\end{equation}
for $n=0,1,2,\dots$.
Furthermore, we can find $C>0$ such that
\begin{equation*}
  \log a_{n+1} \ge pq \log a_n - C(n+1).
\end{equation*}
Then there exists a constant $a_*>0$ such that
\begin{equation}
  \label{eq:3.12}
  a_n \ge a_*^{(pq)^n} M_{z,\tau}(u)^{(pq)^n},
  \qquad n=0,1,2,\dots.
\end{equation}
(See also \cite{HI01} and \cite[Lemma~2.20~(i)]{LN}.)
Combining \eqref{eq:3.10}, \eqref{eq:3.11} and \eqref{eq:3.12},
we see that
\begin{align*}
  \infty>U(x,t) & \ge
  a_*^{(pq)^n} M_{z,\tau}(u)^{(pq)^n}t^{[(pq)^n-1]\frac{p+1}{pq-1}}
  G(x,t+\rho^2)^{(pq)^n}\\
  & =
  \left[a_* M_{z,\tau}(u)t^{\frac{p+1}{pq-1}}
  G(x,t+\rho^2)\right]^{(pq)^n}t^{-\frac{p+1}{pq-1}},
  \quad n=0,1,2,\dots,
\end{align*}
for almost all $x\in{\bf R}^N$, $t\in[0,2\rho^2]$ and $\tau\in(0,\rho^2]$.
This implies that
$$
a_* M_{z,\tau}(u)t^{\frac{p+1}{pq-1}}G(x,t+\rho^2)\le 1
$$
for almost all $x\in{\bf R}^N$, $t\in[0,2\rho^2]$ and $\tau\in(0,\rho^2]$.
Then we deduce that
$$
a_* M_{z,\tau}(u)(\rho^2)^{\frac{p+1}{pq-1}}(8\pi\rho^2)^{-\frac{N}{2}}\le 1,
$$
that is,
$$
\int_{B(z,\rho)} u(x,\tau)\,dx=M_{z,\tau}(u)\le (8\pi)^{\frac{N}{2}}a_*^{-1}\rho^{N-\frac{2(p+1)}{pq-1}}
$$
for almost all $\tau\in(0,\rho^2]$.
This completes the proof.
$\Box$\vspace{5pt}
\newline
We complete the proof of Theorem~\ref{Theorem:1.1}.
\vspace{3pt}
\newline
{\bf Proof of Theorem~\ref{Theorem:1.1}.}
It suffices to prove \eqref{eq:1.9} under the assumption that $T=1$ and $D=1$.
First, we prove that
\begin{equation}
\label{eq:3.13}
\sup_{z\in{\bf R}^N}\mu(B(z,\sigma))\le C\sigma^{N-\frac{2(p+1)}{pq-1}}
\end{equation}
for $0<\sigma\le 1$.
By Lemma~\ref{Lemma:3.2} we have
\begin{equation}
\label{eq:3.14}
\sup_{z\in{\bf R}^N}\int_{B(z,\rho)} u(x,\tau)\,dx\le C_*\rho^{N-\frac{2(p+1)}{pq-1}}
\end{equation}
for all $\rho\in(0,1/\sqrt{5})$ and almost all $\tau\in(0,\rho^2]$,
where $C_*$ is as in Lemma~\ref{Lemma:3.2}.
Let $\zeta\in C_0({\bf R}^N)$ be such that
$$
0\le\zeta\le 1\quad\mbox{in}\quad{\bf R}^N,
\qquad
\zeta=1\quad\mbox{in}\quad B(0,\rho/2),
\qquad
\zeta=0\quad\mbox{outside}\quad B(0,\rho).
$$
Then, by \eqref{eq:1.8} and \eqref{eq:3.14}
we have
\begin{equation*}
\begin{split}
C_*\rho^{N-\frac{2(p-1)}{pq-1}} & \ge \int_{B(z,\rho)} u(x,\tau)\,dx
\ge \int_{{\bf R}^N} u(x,\tau)\zeta(x-z)\,dx\\
 & \to \int_{{\bf R}^N}\zeta(x-z)\,d\mu(x)
\ge\mu(B(z,\rho/2))
\end{split}
\end{equation*}
as $\tau\to 0$, for $z\in{\bf R}^N$ and $\rho\in(0,1/\sqrt{5})$.
This together with Lemma~\ref{Lemma:2.4} implies that
$$
\sup_{z\in{\bf R}^N}\mu(B(z,2\sqrt{5}\rho))
\le C\sup_{z\in{\bf R}^N}\mu(B(z,\rho/2))
\le C\rho^{N-\frac{2(p+1)}{pq-1}}
$$
for $\rho\in(0,1/\sqrt{5})$. Then we set $\sigma:=2\sqrt{5}\rho$ and obtain \eqref{eq:3.13}.

Next we prove that
\begin{equation}
\label{eq:3.15}
\sup_{z\in{\bf R}^N}\nu(B(z,\sigma))\le C\sigma^{N-\frac{2(q+1)}{pq-1}}
\end{equation}
for $0<\sigma\le 1$.
It follows from Lemma~\ref{Lemma:3.2} with $\tau=r^2$ that
\begin{equation}
\label{eq:3.16}
\sup_{z\in{\bf R}^N}\int_{B(z,r)} u(x,r^2)\,dx\le C_*r^{N-\frac{2(p+1)}{pq-1}}
\end{equation}
for almost all $r\in(0,1/\sqrt{5})$.
Let $z\in{\bf R}^N$.
Since Lemma~\ref{Lemma:2.2} implies that
\begin{equation*}
\begin{split}
\int_{{\bf R}^N}G(x-y,t-s)G(y,s+\rho^2)^p\,dy
& \ge C(s+\rho^2)^{-\frac{N}{2}(p-1)}G\left(x,\frac{t+\rho^2}{\max\{p,1\}}\right)
\end{split}
\end{equation*}
for $x\in{\bf R}^N$, $\rho>0$ and $\rho^2<s<t$,
by \eqref{eq:3.4}, \eqref{eq:3.5} and \eqref{eq:3.7}
we have
\begin{equation*}
  \begin{split}
    u(x,\tau+4\rho^2)
    & =U(x-z,2\rho^2)\\
    &
    \ge (D')^{-\frac{N}{2}}\int_0^{2\rho^2}
    \int_{{\bf R}^N}G(x-z-y,2\rho^2-s)V(y,s)^p\,dy\,ds\\
    &
    \ge CM_{z,\tau}(v)^p
    \int_{\rho^2}^{2\rho^2}\int_{{\bf R}^N}G(x-z-y,2\rho^2-s)G(y,s+\rho^2)^p\,dy\,ds\\
    &
    \ge CM_{z,\tau}(v)^p G\left(x-z,\frac{3\rho^2}{\max\{p,1\}}\right)
    \int_{\rho^2}^{2\rho^2} (s+\rho^2)^{-\frac{N}{2}(p-1)}\,ds\\
    &
    \ge CM_{z,\tau}(v)^p \rho^{-N(p-1)+2}G\left(x-z,\frac{3\rho^2}{\max\{p,1\}}\right)
  \end{split}
\end{equation*}
for almost all $x\in{\bf R}^N$, $\rho\in(0,1/5)$ and $\tau\in(0,\rho^2]$.
Since $\tau+4\rho^2\le 5\rho^2<1/5$ for $\rho\in (0,1/5)$,
by \eqref{eq:3.16} we obtain
\begin{equation*}
\begin{split}
C_*(\tau+4\rho^2)^{\frac{N}{2}-\frac{p+1}{pq-1}}
 & \ge\int_{B(z,\sqrt{\tau+4\rho^2})}u(x,\tau+4\rho^2)\,dx
\ge\int_{B(z,2\rho)} u(x,\tau+4\rho^2)\,dx\\
 & \ge CM_{z,\tau}(v)^p \rho^{-N(p-1)+2}
\int_{B(0,\rho)}G\left(x,\frac{3\rho^2}{\max\{p,1\}}\right)\,dx\\
 & \ge CM_{z,\tau}(v)^p \rho^{-N(p-1)+2}
\end{split}
\end{equation*}
for almost all $\rho\in(0,1/5)$ and $\tau\in(0,\rho^2]$.
This implies that
\begin{equation}
  \label{eq:3.17}
  \sup_{z\in{\bf R}^N}\int_{B(z,\rho)}v(x,\tau)\,dx
  = \sup_{z\in{\bf R}^N} M_{z,\tau}(v)
  \le C\rho^{N-\frac{2(q+1)}{pq-1}}
\end{equation}
for all $\rho\in(0,1/5)$ and almost all $\tau\in(0,\rho^2]$.
Then we apply the same argument as in the proof of \eqref{eq:3.13} to obtain \eqref{eq:3.15}.
Therefore, by \eqref{eq:3.13} and \eqref{eq:3.15}
we have \eqref{eq:1.9}, and the proof of Theorem~\ref{Theorem:1.1} is complete.
$\Box$
\section{Proof of Theorem~\ref{Theorem:1.2}}
Applying similar arguments as in the proof of \cite[Theorem~1.2]{HI01},
we prove Theorem~\ref{Theorem:1.2}.
\vspace{3pt}
\newline
{\bf Proof of Theorem~\ref{Theorem:1.2}.}
Similarly to the proof of Theorem~\ref{Theorem:1.1},
it suffices to consider the case of $T=1$.
We prove assertion~(i).
By Theorem~\ref{Theorem:1.1}
we can find a pair of Radon measures $(\mu',\nu')$ satisfying
\begin{equation}
\label{eq:4.1}
\begin{split}
 & \underset{t\to+0}{\mbox{{\rm ess lim}}}\,\int_{{\bf R}^N} u(y,t)\eta(y)\,dy=\int_{{\bf R}^N}\eta(y)\,d\mu'(y),\\
 & \underset{t\to+0}{\mbox{{\rm ess lim}}}\,\int_{{\bf R}^N} v(y,t)\eta(y)\,dy=\int_{{\bf R}^N}\eta(y)\,d\nu'(y),
\end{split}
\end{equation}
for all $\eta\in C_0({\bf R}^N)$.
Furthermore,
for any $r>0$,
by Lemmas~\ref{Lemma:2.4} and \ref{Lemma:3.2} we see that
\begin{equation}
\label{eq:4.2}
\mbox{$\displaystyle{\sup_{z\in {\bf R}^N} \int_{B(z,r)}u(y,\tau) \, dy}$
is uniformly bounded for almost all $\tau\in(0,1/5)$}.
\end{equation}

Let $\zeta\in C_0({\bf R}^N)$ be such that $\zeta\ge 0$ in ${\bf R}^N$.
Let $R>0$ be such that $\mbox{supp}\,\zeta\subset B(0,R)$.
It follows from Definition~\ref{Definition:1.1} and Remark~\ref{Remark:1.1} that
\begin{equation}
\label{eq:4.3}
\begin{split}
\int_{{\bf R}^N} u(x,t)\zeta(x)\,dx &
= \int_{{\bf R}^N}\int_{{\bf R}^N} G(x-y,D_1(t-\tau))u(y,\tau)\zeta(x)\,dx\,dy\\
& +\int_\tau^t\int_{{\bf R}^N}\int_{{\bf R}^N} G(x-y,D_1(t-s))v(y,s)^p\zeta(x)\,dx\,dy\,ds
\end{split}
\end{equation}
for almost all $0<\tau<t<T$.
Set $\zeta(x,t):=[S(D_1t)\zeta](x)$.
Then
\begin{equation}
\label{eq:4.4}
\lim_{t\to+0}\|\zeta(t)-\zeta\|_{L^\infty({\bf R}^N)}=0.
\end{equation}
By \eqref{eq:4.1} and \eqref{eq:4.2} we have
\begin{equation*}
\begin{split}
 & \int_{{\bf R}^N}\int_{{\bf R}^N} G(x-y,D_1(t-\tau))u(y,\tau)\zeta(x)\,dx\,dy
\ge \int_{B(0,R)} u(y,\tau)\zeta(y,t-\tau)\,dy\\
 & \qquad\quad
\ge \int_{B(0,R)} u(y,\tau)\zeta(y)\,dy-C\|\zeta(t-\tau)-\zeta\|_{L^\infty(B(0,R))}\\
 & \qquad\quad
 \to\int_{{\bf R}^N} \zeta(y)\,d\mu'(y)-C\|\zeta(t)-\zeta\|_{L^\infty(B(0,R))}
 \quad\mbox{as}\quad\tau\to +0.
\end{split}
\end{equation*}
This together with \eqref{eq:4.3} implies that
\begin{equation}
\label{eq:4.5}
\begin{split}
\int_{{\bf R}^N} u(x,t)\zeta(x)\,dx & \ge \int_{{\bf R}^N} \zeta(y)\,d\mu'(y)-C\|\zeta(t)-\zeta\|_{L^\infty(B(0,R))}\\
 & +\int_0^t\int_{{\bf R}^N}\int_{{\bf R}^N} G(x-y,D_1(t-s))v(y,s)^p\zeta(x)\,dx\,dy\,ds
\end{split}
\end{equation}
for almost all $t\in(0,T)$.
By \eqref{eq:4.1}, \eqref{eq:4.4} and \eqref{eq:4.5}, letting $t\to+0$,
we obtain
\begin{equation}
\label{eq:4.6}
\underset{t\to+0}{\mbox{{\rm ess lim}}}\,\int_0^t\int_{{\bf R}^N}
\int_{{\bf R}^N} G(x-y,D_1(t-s))v(y,s)^p\zeta(x)\,dx\,dy\,ds=0.
\end{equation}
On the other hand, it follows from Definition~\ref{Definition:1.1}~(ii) that
\begin{equation}
\label{eq:4.7}
\begin{split}
\int_{{\bf R}^N} u(x,t)\zeta(x)\,dx &
=\int_{{\bf R}^N}\int_{{\bf R}^N} G(x-y,D_1t)\zeta(x)\,d\mu(y)\,dx\\
& +\int_0^t\int_{{\bf R}^N}\int_{{\bf R}^N} G(x,y,D_1(t-s))v(y,s)^p\zeta(x)\,dx\,dy\,ds
\end{split}
\end{equation}
for almost all $0<t<T$.
Then, by \eqref{eq:4.1}, \eqref{eq:4.4}, \eqref{eq:4.6} and \eqref{eq:4.7} we obtain
$$
\int_{{\bf R}^N}\zeta(x)\,d\mu'(x)
=\underset{t\to+0}{\mbox{{\rm ess lim}}}\,\int_{{\bf R}^N} \zeta(y,t)\,d\mu(y)
=\int_{{\bf R}^N}\zeta(y)\,d\mu(y).
$$
Since $\zeta$ is arbitrary, we deduce that $\mu=\mu'$ in ${\bf R}^N$.
Similarly, we see that $\nu=\nu'$ in ${\bf R}^N$.
Then assertion~(i) follows.

We prove assertion~(ii).
Let $(u,v)$ be a solution to \eqref{eq:1.1} in ${\bf R}^N\times(0,1)$.
Let $0<t<1$ and $n=1,2,\dots$.
By the Besicovitch covering lemma we can find an integer $m$
depending only on $N$ and a  set $\{x_{k,i}\}_{k=1,\dots,m,\,i\in{\bf N}}\subset{\bf R}^N\setminus B(0,nt^{1/2})$ such that
\begin{equation}
\label{eq:4.8}
B_{k,i}\cap B_{k,j}=\emptyset\quad\mbox{if $i\not=j$}
\qquad\mbox{and}\qquad
{\bf R}^N\setminus B(0,nt^{\frac{1}{2}})\subset\bigcup_{k=1}^m\bigcup_{i=1}^\infty B_{k,i},
\end{equation}
where $B_{k,i}:=\overline{B(x_{k,i},t^{1/2})}$.
By Lemma~\ref{Lemma:2.4}, Lemma~\ref{Lemma:3.2} and \eqref{eq:4.8} we obtain
\begin{equation}
\label{eq:4.9}
\begin{split}
 & \underset{0<\tau<t/5}{\mbox{{\rm ess sup}}}
 \int_{{\bf R}^N\setminus B(0,nt^\frac{1}{2})}G(y,D_1(t-\tau))u(y,\tau)\,dy\\
 &  \le\sum^m_{k=1}\sum^\infty_{i=1}\,\underset{0<\tau<t/5}{\mbox{{\rm ess sup}}}\int_{B_{k,i}}G(y,D_1(t-\tau))u(y,\tau)\,dy\\
 & \le C\,\underset{0<\tau<t/5}{\mbox{{\rm ess sup}}}\sup_{z\in{\bf R}^N}\int_{B(z,t^\frac{1}{2})}u(y,\tau)\,dy\\
 & \qquad\qquad
 \times\sum^m_{k=1}\sum^\infty_{i=1}\sup_{0<\tau<t/5}\sup_{y\in B_{k,i}}
 (t-\tau)^{-\frac{N}{2}}\exp\left(-\frac{|y|^2}{4D_1(t-\tau)}\right)\\
  & \le C t^{-\frac{p+1}{pq-1}} \sum^m_{k=1}\sum^\infty_{i=1}\sup_{y\in B_{k,i}}
 \exp\left(-\frac{|y|^2}{4D_1t}\right).
\end{split}
\end{equation}
On the other hand,
since
\begin{equation*}
(n-1)t^{\frac{1}{2}}\le |z| \le |z-x_{k,i}| + |x_{k,i}-y| + |y| \le |y| + 2t^\frac{1}{2}
\end{equation*}
for $y$, $z \in B_{k,i}$,
we have
$$
\inf_{y\in B_{k,i}}\frac{|y|^2}{t}
\ge \frac{||z|-2t^\frac{1}{2}|^2}{t}
\ge\frac{|z|^2}{2t}-C
$$
for $z\in B_{k,i}$ and sufficiently large $n$.
Then
$$
\sup_{y\in B_{k,i}}
\exp\left(-\frac{|y|^2}{4D_1t}\right)
\le C\,\, \frac{1}{|B_{k,i}|}\int_{B_{k,i}}
\exp\left(-\frac{|z|^2}{8D_1t}\right)\,dz
$$
for sufficiently large $n$,
where $|B_{k,i}|$ denotes the measure of the ball $B_{k,i}$.
This together with \eqref{eq:4.8} and \eqref{eq:4.9} implies that
\begin{equation}
\label{eq:4.10}
\begin{split}
 & \underset{0<\tau<t/5}{\mbox{{\rm ess sup}}}\,
 \int_{{\bf R}^N\setminus B(0,nt^\frac{1}{2})}G(y,D_1(t-\tau))u(y,\tau)\,dy\\
 & \le C t^{-\frac{N}{2}-\frac{p+1}{pq-1}}
 \sum^m_{k=1}\sum^\infty_{i=1}\int_{B_{k,i}}\exp\left(-\frac{|z|^2}{8D_1t}\right)\,dz\\
 & \le C t^{-\frac{N}{2}-\frac{p+1}{pq-1}}
 \int_{{\bf R}^N\setminus B(0,(n-1)t^{\frac{1}{2}})}\exp\left(-\frac{|z|^2}{8D_1t}\right)\,dz\\
 & \le C t^{-\frac{p+1}{pq-1}}
 \int_{{\bf R}^N\setminus B(0,n-1)}\exp\left(-\frac{|z|^2}{8D_1}\right)\,dz\to 0
\end{split}
\end{equation}
as $n\to\infty$.
Similarly, by Theorem~\ref{Theorem:1.1} we have
\begin{equation}
\label{eq:4.11}
\begin{split}
 & \int_{{\bf R}^N\setminus B(0,nt^\frac{1}{2})}G(y,D_1t)\,d\mu(y) \\
 & \le C\sup_{z\in{\bf R}^N}\mu(B(z,t^\frac{1}{2}))
 \sum^m_{k=1}\sum^\infty_{i=1}  \sup_{y\in B_{k,i}}
 t^{-\frac{N}{2}}\exp\left(-\frac{|y|^2}{4D_1t}\right)\\
 & \le C t^{-\frac{p+1}{pq-1}}
 \int_{{\bf R}^N\setminus B(0,n-1)}\exp\left(-\frac{|z|^2}{8D_1}\right)\,dz\to 0
\end{split}
\end{equation}
as $n\to\infty$.
In particular,
by \eqref{eq:4.2}, \eqref{eq:4.10} and \eqref{eq:4.11}
we see that
\begin{equation}
\label{eq:4.12}
\int_{{\bf R}^N}G(y,D_1(t-\tau))u(y,\tau)\,dy<\infty,
\qquad
\int_{{\bf R}^N} G(y,D_1t)\,d\mu(y)<\infty,
\end{equation}
for all $0<t<1$ and almost all $\tau\in(0,t/5)$.

Let $\eta_n\in C_0({\bf R}^N)$ be such that
$$
0\le\eta_n\le 1\quad\mbox{in}\quad{\bf R}^N,
\qquad
\eta_n=1\quad\mbox{on}\quad B(0,nt^\frac{1}{2}),
\qquad
\eta_n=0\quad\mbox{outside}\quad B(0,2nt^\frac{1}{2}).
$$
It follows from \eqref{eq:4.12} that
\begin{equation}
\label{eq:4.13}
\begin{split}
 & \left|\int_{{\bf R}^N}G(y,D_1(t-\tau))u(y,\tau)\,dy-\int_{{\bf R}^N} G(y,D_1t)\,d\mu(y)\right|\\
 & \le\left|\int_{{\bf R}^N}G(y,D_1t)u(y,\tau)\eta_n(y)\,dy - \int_{{\bf R}^N} G(y,D_1t)\eta_n(y)\,d\mu(y)\right|\\
 & \qquad
 +\left|\int_{{\bf R}^N}[G(y,D_1(t-\tau))-G(y,D_1t)]u(y,\tau)\eta_n(y)\,dy\right|\\
 & \qquad\qquad
 +\int_{{\bf R}^N\setminus B(0,nt^\frac{1}{2})}G(y,D_1(t-\tau))u(y,\tau)\,dy
 +\int_{{\bf R}^N\setminus B(0,nt^\frac{1}{2})} G(y,D_1t)\,d\mu(y)
\end{split}
\end{equation}
for $n=1,2,\dots$ and almost all $\tau\in(0,t/2)$.
By \eqref{eq:1.8}
we see that
\begin{equation}
  \label{eq:4.14}
  \underset{\tau\to+0}{\mbox{{\rm ess lim}}}\,
  \left[\int_{{\bf R}^N}G(y,D_1t)u(y,\tau)\eta_n(y)\,dy - \int_{{\bf R}^N} G(y,D_1t)\eta_n(y)\,d\mu(y)\right]=0.
\end{equation}
Furthermore, by \eqref{eq:4.2} we have
\begin{equation}
  \label{eq:4.15}
  \begin{split}
    & \underset{\tau\to+0}{\mbox{{\rm ess lim}}}\,
    \left|\int_{{\bf R}^N}[G(y,D_1(t-\tau))-G(y,D_1t)]u(y,\tau)\eta_n(y)\,dy\right|\\
    &
    \le \sup_{y\in B(0,2nt^{\frac{1}{2}}),s\in(4t/5,t)}
    \,|\partial_t G(y,D_1s)|\,\underset{\tau\to+0}{\mbox{{\rm ess limsup}}}\,
    \biggr[\tau\int_{B(0,2nt^\frac{1}{2})}u(y,\tau)\,dy\biggr]=0.
  \end{split}
\end{equation}
By \eqref{eq:4.13}, \eqref{eq:4.14} and \eqref{eq:4.15} we obtain
\begin{equation*}
  \begin{split}
    & \underset{\tau\to+0}{\mbox{{\rm ess lim}}}\,
    \left|\int_{{\bf R}^N}G(y,D_1(t-\tau))u(y,\tau)\,dy - \int_{{\bf R}^N} G(y,D_1t)\,d\mu(y)\right|\\
    & \le
    \underset{0<\tau<t/5}{\mbox{{\rm ess sup}}}\int_{{\bf R}^N\setminus B(0,nt^\frac{1}{2})}G(y,D_1(t-\tau))u(y,\tau)\,dy
    +\int_{{\bf R}^N\setminus B(0,nt^\frac{1}{2})} G(y,D_1t)\,d\mu(y)
  \end{split}
\end{equation*}
for $n=1,2,\dots$.
This together with \eqref{eq:4.10} and \eqref{eq:4.11} implies that
$$
\underset{\tau\to+0}{\mbox{{\rm ess lim}}}\,
\left|\int_{{\bf R}^N}G(y,D_1(t-\tau))u(y,\tau)\,dy - \int_{{\bf R}^N} G(y,D_1t)\,d\mu(y)\right|=0.
$$
Similarly, we see that
$$
\underset{\tau\to+0}{\mbox{{\rm ess lim}}}\,
\left|\int_{{\bf R}^N}G(y,D_2(t-\tau))v(y,\tau)\,dy - \int_{{\bf R}^N} G(y,D_2t)\,d\nu(y)\right|=0.
$$
These together with Definition~\ref{Definition:1.1}~(i) imply that
$(u,v)$ is a solution to \eqref{eq:1.1} with \eqref{eq:1.2} in ${\bf R}^N\times(0,1)$.
Thus assertion~(ii) follows, and the proof of Theorem~\ref{Theorem:1.2} is complete.
$\Box$
\section{Initial traces in case (B)}
In this section we focus on case~(B), that is,
\begin{equation}
\label{eq:5.1}
0<p<q,\qquad pq>1,\qquad
\frac{q+1}{pq-1}=\frac{N}{2},
\end{equation}
and prove assertion~(1) of Theorem~\ref{Theorem:1.3}.
It follows from \eqref{eq:5.1} that
\begin{equation}
\label{eq:5.2}
-\frac{N}{2}(p-1)+1>0>-\frac{N}{2}(q-1)+1.
\end{equation}
Similarly to Section~2,
it suffices to consider the case where $T=1$ and $D=1$.
\vspace{3pt}

Let
\begin{equation}
\label{eq:5.3}
z\in{\bf R}^N,\qquad
0<\rho<1/\sqrt{10},\qquad
0<\tau\le\rho^2<1.
\end{equation}
Similarly to \eqref{eq:3.4},
we define $U=U(x,t)$ and $V=V(x,t)$ for $x\in{\bf R}^N$ and $0<t<t_*$, where
\begin{equation}
\label{eq:5.4}
t_*:=1-\tau-2\rho^2>\frac{1}{2}>2\rho^2.
\end{equation}
Then $U$ and $V$ satisfy \eqref{eq:3.5} and \eqref{eq:3.7}
for almost all $x\in{\bf R}^N$, $t\in(0,t_*)$ and $\tau\in(0,\rho^2]$.
Furthermore, the following lemma holds.
\begin{lemma}
\label{Lemma:5.1}
Assume \eqref{eq:5.1} and \eqref{eq:5.3}.
Let $(U,V)$ satisfy \eqref{eq:3.5}.
Assume that there exist constants $a>0$, $b\ge 1$ and $c\ge 0$ such that
\begin{equation}
\label{eq:5.5}
V(x,t) \ge aG\left(x,\frac{t+\rho^2}{b}\right)\left[\log\frac{t}{\rho^2}\right]^{c}
\end{equation}
for almost all $x\in{\bf R}^N$ and $t\in[2\rho^2,t_*)$,
where $\rho>0$.
Then there exists a positive constant $\gamma_2$ depending only on $N$, $D_1$, $D_2$, $p$ and $q$ such that
$$
V(x,t) \ge \gamma_2
a^{pq}b^{\frac{N}{2}q(p-2)+\frac{N}{2}(q-2)} (pc+1)^{-q}(pqc+1)^{-1}
G\left(x,\frac{t+\rho^2}{q\max\{pb,1\}}\right)\left[\log\frac{t}{\rho^2}\right]^{pqc+1}
$$
for almost all $x\in{\bf R}^N$ and $t\in[2\rho^2,t_*)$.
\end{lemma}
For the proof of Lemma~\ref{Lemma:5.1},
we prove:
\begin{Lemma}
\label{Lemma:5.2}
Let $a>-1$, $b>0$ and $\rho>0$.
Then
$$
\int_{\rho^2}^t (s+\rho^2)^a \left[ \log \frac{s}{\rho^2} \right]^b \, ds
\ge \frac{(t+\rho^2)^{a+1}}{4(a+b+2)}\left[ \log \frac{t}{\rho^2} \right]^b
\quad\mbox{for $t\ge {2\rho^2}$}.
$$
\end{Lemma}
{\bf Proof.}
Let $t\ge 2\rho^2$.
Set
\begin{equation*}
  F(t)
  :=
  \int_{\rho^2}^t (s+\rho^2)^a \left[ \log \frac{s}{\rho^2} \right]^b \, ds
  - \frac{1}{4(a+b+2)} (t+\rho^2)^{a+1} \left[ \log \frac{t}{\rho^2} \right]^b.
\end{equation*}
Then, since
\begin{align*}
  \int_1^2 (\sigma+1)^a (\log \sigma)^b \, d\sigma
  & =
  \int_1^2 (\sigma+1)^{a+1} \cdot \frac{\sigma}{\sigma+1} \cdot \sigma^{-1} (\log \sigma)^b \, d\sigma
  \\
  & \ge
  \frac{1}{2} \cdot 2^{a+1} \int_1^2 \sigma^{-1} (\log \sigma)^b \, d\sigma
  =
  \frac{1}{2} \cdot 2^{a+1} \cdot \frac{1}{b+1} (\log 2)^{b+1}
  \\
  & \ge
  (2\log 2) \cdot \frac{1}{4(a+b+2)} 2^{a+1} (\log 2)^{b}
  \ge \frac{1}{4(a+b+2)} 2^{a+1} (\log 2)^{b},
\end{align*}
we have
\begin{equation}
  \label{eq:a}
  \begin{aligned}
    F(2\rho^2)
    & =
    \int_{\rho^2}^{2\rho^2} (s+\rho^2)^a \left[ \log \frac{s}{\rho^2} \right]^b \, ds
    - \frac{1}{4(a+b+2)} (2\rho^2)^{a+1} (\log 2)^b
    \\
    & = (\rho^2)^{a+1} \left\{
    \int_1^2 (\sigma+1)^a (\log \sigma)^b \, d\sigma
    - \frac{1}{4(a+b+2)} 2^{a+1} (\log 2)^{b}
    \right\}
    \ge 0.
  \end{aligned}
\end{equation}
Since
\begin{align*}
  F'(t)
  & =
  (t+\rho^2)^{a} \left[ \log \frac{t}{\rho^2} \right]^b
  \\
  & \qquad
  - \frac{1}{4(a+b+2)} \left\{
  (a+1) (t+\rho^2)^{a} \left[ \log \frac{t}{\rho^2} \right]^b
  + b (t+\rho^2)^{a}\cdot \frac{t+\rho^2}{t} \left[ \log \frac{t}{\rho^2} \right]^{b-1}
  \right\}
  \\
  & =
  \left\{
  1 - \frac{a+1}{4(a+b+2)} - \frac{b}{4(a+b+2)}
  \cdot \frac{t+\rho^2}{t} \left[ \log \frac{t}{\rho^2} \right]^{-1}
  \right\}
  (t+\rho^2)^{a} \left[ \log \frac{t}{\rho^2} \right]^b
  \\
  & \ge
  \left\{ \frac{3}{4}
  - \frac{1}{4} \cdot \frac{3}{2} ( \log 2 )^{-1} \right\}
  (t+\rho^2)^{a} \left[ \log \frac{t}{\rho^2} \right]^b
  \ge 0
\end{align*}
for $t\ge 2\rho^2$,
by \eqref{eq:a} we have $F(t) \ge 0$ for $t\ge 2\rho^2$.
Thus Lemma~\ref{Lemma:5.2} follows.
$\Box$\vspace{5pt}
\newline
{\bf Proof of Lemma~\ref{Lemma:5.1}.}
By Lemma~\ref{Lemma:2.2}, Lemma~\ref{Lemma:5.2}, \eqref{eq:3.5}, \eqref{eq:5.2} and \eqref{eq:5.5} we have
\begin{align*}
  U(x,t)
  &
  \ge (D')^{-\frac{N}{2}}\int_{\rho^2}^t\int_{{\bf R}^N}G(x-y,t-s)V(y,s)^p\,dy\,ds
  \\
  &
  \ge Ca^p\int_{\rho^2}^t\int_{{\bf R}^N}
  \left[\log\frac{s}{\rho^2}\right]^{pc}
  G(x-y,t-s)G\left(y,\frac{s+\rho^2}{b}\right)^p\,dy\,ds
  \\
  &
  \ge Ca^pb^{\frac{N}{2}(p-1)}
  \left(\frac{\min\{pb,1\}}{\max\{pb,1\}}\right)^{\frac{N}{2}}G\left(x,\frac{t+\rho^2}{\max\{pb,1\}}\right)
  \\
  & \hspace{4.5 true cm}
  \times \int_{\rho^2}^t
  (s+\rho^2)^{-\frac{N}{2}(p-1)}\left[\log\frac{s}{\rho^2}\right]^{pc}\,ds
  \\
  & \ge C\frac{a^pb^{\frac{N}{2}(p-2)}}{pc+1} (t+\rho^2)^{-\frac{N}{2}(p-1)+1}
  G\left(x,\frac{t+\rho^2}{\max\{pb,1\}}\right)
  \left[\log\frac{t}{\rho^2}\right]^{pc}
\end{align*}
for almost all $x\in{\bf R}^N$ and $t\in[{2\rho^2},t_*)$.
Since
\[
-\frac{N}{2}(pq-q)+q -\frac{N}{2}(q-1)
= -\frac{N}{2}(pq-1) + q
= -\frac{q+1}{pq-1}(pq-1) + q
= -1,
\]
by Lemma~\ref{Lemma:2.2} we obtain
\begin{equation*}
  \begin{aligned}
    V(x,t)
    & \ge (D')^{-\frac{N}{2}}\int_{\rho^2}^t \int_{{\bf R}^N}G(x-y,t-s)U(y,s)^q\,dy\,ds
    \\
    & \ge (D')^{-\frac{N}{2}}\left(C\frac{a^pb^{\frac{N}{2}(p-2)}}{pc+1}\right)^q
    \int_{\rho^2}^t (s+\rho^2)^{-\frac{N}{2}(pq-q)+q}\left[\log\frac{s}{\rho^2}\right]^{pqc}
    \\
    & \hspace{3.0 true cm}
    \times
    \left[
    \int_{{\bf R}^N} G(x-y,t-s)G\left(x,\frac{s+\rho^2}{\max\{pb,1\}}\right)^q dy
    \right] ds
    \\
    & \ge C\left(\frac{a^pb^{\frac{N}{2}(p-2)}}{pc+1}\right)^q
    [\max\{pb,1\}]^{\frac{N}{2}(q-2)}
    G\left(x,\frac{t+\rho^2}{q\max\{pb,1\}}\right)
    \\
    & \hspace{3.0 true cm}
    \times \int_{\rho^2}^t (s+\rho^2)^{-\frac{N}{2}(pq-q)+q -\frac{N}{2}(q-1)}\left[\log\frac{s}{\rho^2}\right]^{pqc} ds
    \\
    & \ge C\frac{a^{pq} b^{\frac{N}{2}q(p-2)+\frac{N}{2}(q-2)}}{(pc+1)^q}
    G\left(x,\frac{t+\rho^2}{q\max\{pb,1\}}\right)
    \int_{\rho^2}^t (s+\rho^2)^{-1}\left[\log\frac{s}{\rho^2}\right]^{pqc} ds
    \\
    & \ge C\frac{a^{pq}b^{\frac{N}{2}q(p-2)+\frac{N}{2}(q-2)}}{(pc+1)^q(pqc+1)}
    G\left(x,\frac{t+\rho^2}{q\max\{pb,1\}}\right)\left[\log\frac{t}{\rho^2}\right]^{pqc+1}
  \end{aligned}
\end{equation*}
for almost all $x\in{\bf R}^N$ and $t\in[{2\rho^2},t_*)$.
Thus Lemma~\ref{Lemma:5.1} follows.
$\Box$\vspace{5pt}
\newline
Combining Lemma~\ref{Lemma:5.1} and \eqref{eq:3.7}, we have:
\begin{lemma}
\label{Lemma:5.3}
Assume \eqref{eq:5.1}.
Let $(u,v)$ satisfy \eqref{eq:3.2}.
Then there exists a constant $C_*>0$ such that
\begin{equation}
\label{eq:5.6}
\sup_{z\in{\bf R}^N}\int_{B(z,\rho)}v(x,\tau)\, dx \le C_*\left[\log\biggr(e+\frac{1}{\rho}\biggr)\right]^{-\frac{1}{pq-1}}
\end{equation}
for all $\rho\in(0,1/\sqrt{10})$ and almost all $\tau\in(0,\rho^2]$.
\end{lemma}
{\bf Proof.}
Let $z\in \mathbf{R}^N$ and $\rho\in(0,1/\sqrt{10})$.
Combining Lemma~\ref{Lemma:5.1} and \eqref{eq:3.7},
we obtain
\begin{equation}
\label{eq:5.7}
\infty>V(x,t)\ge a_nG\left(x,\frac{t+\rho^2}{b_n}\right)\left[\log\frac{t}{\rho^2}\right]^{c_n},
\qquad n=0,1,2,\dots,
\end{equation}
for almost all $x\in{\bf R}^N$, $t\in[{2\rho^2},t_*)$ and $\tau\in(0,\rho^2]$.
Here $\{a_n\}_{n=0}^\infty$, $\{b_n\}_{n=0}^\infty$ and $\{c_n\}_{n=0}^\infty$ are sequences defined by
$$
a_{n+1}:=\gamma_2
\frac{a_n^{pq}b_n^{\frac{N}{2}q(p-2)+\frac{N}{2}(q-2)}}{(pc_n+1)^q(pqc_n+1)},\quad
b_{n+1}:=q\max\{pb_n,1\},\quad
c_{n+1}:=pq c_n+1,
$$
for $n=0,1,2,\dots$ with
$$
a_0:=c_*M_{z,\tau}(v),
\qquad
b_0=1,
\qquad
c_0=0,
$$
where $\gamma_2$ is as in Lemma~\ref{Lemma:5.1}.
Since $pb_1=pq\max\{p,1\}>1$, we have
\begin{equation}
\label{eq:5.8}
b_n=(pq)^{n-1}b_1,\qquad c_n=\frac{(pq)^n-1}{pq-1},
\end{equation}
for $n=1,2,3,\dots$.
Furthermore, similarly to \eqref{eq:3.12},
we can find $a_*>0$ such that
\begin{equation}
\label{eq:5.9}
a_n\ge a_*^{(pq)^n}M_{z,\tau} (v)^{(pq)^n},\qquad n=1,2,\dots.
\end{equation}
By \eqref{eq:5.7}, \eqref{eq:5.8} and \eqref{eq:5.9}
we see that
\begin{equation*}
\begin{split}
\infty>V(x,t) & \ge a_*^{(pq)^n}M_{z,\tau} (v)^{(pq)^n}
G\left(x,\frac{t+\rho^2}{(pq)^{n-1}b_1}\right)\left[\log\frac{t}{\rho^2}\right]^{\frac{(pq)^n-1}{pq-1}}\\
 & \ge\left\{a_*M_{z,\tau}(v)\left[\log\frac{t}{\rho^2}\right]^{\frac{1}{pq-1}}\right\}^{(pq)^n}
G\left(x,\frac{t+\rho^2}{(pq)^{n-1}b_1}\right)\left[\log\frac{t}{\rho^2}\right]^{\frac{-1}{pq-1}}
\end{split}
\end{equation*}
for almost all $x\in{\bf R}^N$, $t\in[{2\rho^2},t_*)$ and $\tau\in(0,\rho^2]$
and for all $n=1,2,\dots$.
This implies that
$$
a_*M_{z,\tau}(v)\left[\log\frac{t}{\rho^2}\right]^{\frac{1}{pq-1}}\le 1
$$
for almost all $t\in[{2\rho^2},t_*)$ and $\tau\in(0,\rho^2]$.
Since $t_*>1/2$, we deduce that
$$
\int_{B(z,\rho)}v(x,\tau)\,dx=M_{z,\tau}(v)
\le a_*^{-1}\left[\log\frac{1}{2\rho^2}\right]^{-\frac{1}{pq-1}}
\le C\left[\log\biggr(e+\frac{1}{\rho}\biggr)\right]^{-\frac{1}{pq-1}}
$$
for almost all $\tau\in(0,\rho^2]$.
This implies \eqref{eq:5.6}, and the proof is complete.
$\Box$
\vspace{5pt}

\noindent
\textbf{Proof of assertion~(1) of Theorem~\ref{Theorem:1.3}}.
Applying the same argument as in the proof of Theorem~\ref{Theorem:1.1},
by Lemma~\ref{Lemma:5.3} we observe that
$$
\sup_{z\in{\bf R}^N}\nu(B(z,\sigma))\le C\left[\log\biggr(e+\frac{1}{\sigma}\biggr)\right]^{-\frac{1}{pq-1}}
$$
for $0<\sigma\le 1$.
It remains to prove
\begin{equation}
  \label{eq:5.10}
  \sup_{z\in \mathbf{R}^N} \int_0^\sigma
  \tau^{-1}\left[ \frac{\mu(B(z,\tau))}{\tau^{N-\frac{2(p+1)}{pq-1}}} \right]^q d\tau
  \le C\left[\log\biggr(e+\frac{1}{\sigma}\biggr)\right]^{-\frac{1}{pq-1}}
\end{equation}
for $0<\sigma\le 1$.

Let $\sigma_*=1/\sqrt{10}$ and $0<\sigma<\sigma_*$.
Applying Lemma~\ref{Lemma:5.3} with $\tau=\sigma^2<1/10$,
we have
$$
\sup_{z\in{\bf R}^N}\int_{B(z,\sigma)} v(x,\sigma^2)\, dx \le C\left[\log\biggr(e+\frac{1}{\sigma}\biggr)\right]^{-\frac{1}{pq-1}}
$$
for almost all $\sigma\in(0,\sigma_*)$.
This together with Lemma~\ref{Lemma:2.4} implies that
\begin{equation}
\label{eq:5.11}
\sup_{z\in{\bf R}^N}\int_{B(z,\sqrt{2}\sigma)} v(x,\sigma^2)\, dx \le C\left[\log\biggr(e+\frac{1}{\sigma}\biggr)\right]^{-\frac{1}{pq-1}}
\end{equation}
for almost all $\sigma\in(0,\sigma_*)$.
On the other hand,
by Theorems~\ref{Theorem:1.1} and \ref{Theorem:1.2}
we see that $(u,v)$ is a solution to problem~\eqref{eq:1.1} with \eqref{eq:1.2} in ${\bf R}^N\times(0,1)$.
Then it follows from Definition~\ref{Definition:1.1}~(ii) and \eqref{eq:3.1} that
\begin{equation}
\label{eq:5.12}
\begin{split}
u(x,t) & \ge\int_{{\bf R}^N}G(x-y,D_1t)\,d\mu(y)\\
 & \ge (D')^{-\frac{N}{2}}\int_{{\bf R}^N}G(x-y,t)\,d\mu(y)
=(D')^{-\frac{N}{2}}[S(t)\mu](x)
\end{split}
\end{equation}
for almost all $x\in{\bf R}^N$ and $t\in(0,1)$.
Using Definition~\ref{Definition:1.1}~(ii) again,
by \eqref{eq:3.1} and \eqref{eq:5.12} we have
\begin{equation*}
\begin{split}
v(x,t) & \ge\int_0^t\int_{{\bf R}^N}G(x-y,D_2(t-s))u(y,s)^q\,dy\,ds\\
 & \ge (D')^{-\frac{N}{2}(q+1)}\int_0^t\int_{{\bf R}^N}G(x-y,t-s)[S(s)\mu](y)^q\,dy\,ds
\end{split}
\end{equation*}
for almost all $x\in{\bf R}^N$ and $t\in(0,1)$.
Then we see that
\begin{equation}
\label{eq:5.13}
\begin{split}
 &  \int_{B(z,\sqrt{2}\sigma)} v(x,\sigma^2) \, dx\\
 &  \ge C\int_{B(z,\sqrt{2}\sigma)}
  \int_0^{\sigma^2} \int_{\mathbf{R}^N} G(x-y, \sigma^2-s) [S(s)\mu](y)^q \, dy \, ds\,dx\\
 &  \ge C\int_{B(z,\sqrt{2}\sigma)}
  \int_0^{\sigma^2} \int_{B(z,\sqrt{s})} G(x-y, \sigma^2-s) [S(s)\mu](y)^q \, dy \, ds\,dx\\
 & =C\int_0^{\sigma^2}\int_{B(z,\sqrt{s})}\int_{B(z,\sqrt{2}\sigma)}
   G(x-y, \sigma^2-s) [S(s)\mu](y)^q \,dx\,dy\,ds
\end{split}
\end{equation}
for all $z\in \mathbf{R}^N$ and almost all $\sigma\in(0,\sigma_*)$.
On the other hand,
since $\sqrt{s} + \sqrt{\sigma^2 - s}\le \sqrt{2}\sigma$ for $s\in (0,\sigma^2)$,
it follows that
$$
B(y,\sqrt{\sigma^2-s})\subset B(z,\sqrt{2}\sigma)\quad\mbox{for $y\in B(z,\sqrt{s})$ and $s\in(0,\sigma^2)$}.
$$
This together with \eqref{eq:5.13} implies that
\begin{equation}
\label{eq:5.14}
\begin{split}
 & \int_{B(z,\sqrt{2}\sigma)} v(x,\sigma^2) \, dx\\
 & \ge\int_0^{\sigma^2} \int_{B(z,\sqrt{s})} \int_{B(y,\sqrt{\sigma^2-s})}
    G(x-y,\sigma^2-s)[S(s)\mu](y)^q\,dx \,dy \,ds\\
 & \ge C\int_0^{\sigma^2} \int_{B(z,\sqrt{s})} [S(s)\mu](y)^q \, dy \, ds
\end{split}
\end{equation}
for all $z\in \mathbf{R}^N$ and almost all $\sigma\in(0,\sigma_*)$.
Furthermore,
by Lemma~\ref{Lemma:2.1} we have
\begin{equation}
\label{eq:5.15}
[S(s)\mu](y)\ge \int_{B(z,\sqrt{s})}G(y-\xi,s)\,d\mu(\xi)
\ge CG\left( y-z,\frac{s}{2} \right) \mu(B(z,\sqrt{s}))
\end{equation}
for $y$, $z\in \mathbf{R}^N$ and $s>0$.
Combining \eqref{eq:5.11}, \eqref{eq:5.14} and \eqref{eq:5.15},
we obtain
\begin{equation}
\label{eq:5.16}
\begin{split}
\left[\log\biggr(e+\frac{1}{\sigma}\biggr)\right]^{-\frac{1}{pq-1}}
 & \ge C\int_0^{\sigma^2} \mu(B(z,\sqrt{s}))^q
\int_{B(z,\sqrt{s})} G\left( y-z,\frac{s}{2} \right)^q \,dy\,ds\\
 & \ge C\int_0^{\sigma^2} \mu(B(z,\sqrt{s}))^q s^{-\frac{N}{2}(q-1)}\, ds\\
 & \ge C\int_0^{\sigma} \mu(B(z,\tau))^q \tau^{-Nq+N+1}\, d\tau
\end{split}
\end{equation}
for all $z\in{\bf R}^N$ and $\sigma\in(0,\sigma_*)$.
By \eqref{eq:5.1} we see that
\begin{equation*}
\begin{split}
  & (-Nq+N+1)-\left\{ -q\left(N-\frac{2(p+1)}{pq-1}\right) -1 \right\}\\
  & = N+2-2\cdot\frac{pq+q}{pq-1}
  = N+2-2\cdot \frac{(pq-1)+(q+1)}{pq-1}
  = N - 2\cdot \frac{q+1}{pq-1}
  = 0,
\end{split}
\end{equation*}
which together with \eqref{eq:5.16} implies that
\begin{equation}
\label{eq:5.17}
 \sup_{z\in \mathbf{R}^N} \int_0^\sigma
 \tau^{-1}\left[ \frac{\mu(B(z,\tau))}{\tau^{N-\frac{2(p+1)}{pq-1}}} \right]^q d\tau
 \le C\left[\log\biggr(e+\frac{1}{\sigma}\biggr)\right]^{-\frac{1}{pq-1}}
\end{equation}
for $0<\sigma<\sigma_*=1/\sqrt{10}$.
In the case of $\sigma_*\le\sigma<1$,
by \eqref{eq:5.17} we see that
\begin{equation}
\label{eq:5.18}
\begin{split}
\sup_{z\in \mathbf{R}^N} \int_0^\sigma
\tau^{-1}\left[ \frac{\mu(B(z,\tau))}{\tau^{N-\frac{2(p+1)}{pq-1}}} \right]^q d\tau
 & \le C+\sup_{z\in \mathbf{R}^N} \int_{\sigma_*}^\sigma
 \tau^{-1}\left[ \frac{\mu(B(z,\tau))}{\tau^{N-\frac{2(p+1)}{pq-1}}} \right]^q d\tau\\
 & \le C+C\sup_{z\in \mathbf{R}^N} \int_{\sigma_*/\sqrt{10}}^{\sigma/\sqrt{10}}
 \tau^{-1} \left[ \frac{\mu(B(z,\sqrt{10}\tau))}{\tau^{N-\frac{2(p+1)}{pq-1}}} \right]^q d\tau.
\end{split}
\end{equation}
On the other hand,
by Lemma~\ref{Lemma:2.3}, for any $z\in{\bf R}^N$,
we can find $\{z_i\}_{i=1}^m\subset{\bf R}^N$ such that
$$
\mu(B(z,\sqrt{10}\tau))\le\sum_{i=1}^m\mu(B(z_i,\tau))\quad\mbox{for $\tau\in(0,1)$}.
$$
Here $m$ is independent of $z$ and $\tau$.
Then
$$
\mu(B(z,\sqrt{10}\tau))^q\le C\sum_{i=1}^m \mu(B(z_i,\tau))^q\quad\mbox{for $\tau\in(0,1)$},
$$
which together with \eqref{eq:5.17} implies that
\begin{equation}
\label{eq:5.19}
\int_{\sigma_*/\sqrt{10}}^{\sigma/\sqrt{10}}
\tau^{-1}\left[ \frac{\mu(B(z,\sqrt{10}\tau))}{\tau^{N-\frac{2(p+1)}{pq-1}}} \right]^q d\tau
\le C\sum_{i=1}^m\int_0^{\sigma/\sqrt{10}}
\tau^{-1}\left[ \frac{\mu(B(z_i,\tau))}{\tau^{N-\frac{2(p+1)}{pq-1}}} \right]^q d\tau\le C
\end{equation}
for all $z\in{\bf R}^N$.
We deduce from \eqref{eq:5.18} and \eqref{eq:5.19} that
\begin{equation}
\label{eq:5.20}
\sup_{z\in \mathbf{R}^N} \int_0^\sigma
\tau^{-1}\left[ \frac{\mu(B(z,\tau))}{\tau^{N-\frac{2(p+1)}{pq-1}}} \right]^q d\tau
\le C\le C\left[\log\biggr(e+\frac{1}{\sigma}\biggr)\right]^{-\frac{1}{pq-1}}
\end{equation}
for $\sigma_*\le\sigma<1$.
Combining \eqref{eq:5.17} and \eqref{eq:5.20},
we obtain \eqref{eq:5.10}.
Thus assertion~(1) of Theorem~\ref{Theorem:1.3} follows.
$\Box$
\section{Initial traces in case (C)}
In this section we focus on case (C), that is,
\begin{equation}
\label{eq:6.1}
\frac{q+1}{pq-1} = \frac{N}{2},
\qquad
p = q,
\end{equation}
and prove assertion~(2) of Theorem~\ref{Theorem:1.3}.
It follows from \eqref{eq:6.1} that
\begin{equation}
\label{eq:6.2}
  p = q = 1 + \frac{2}{N}.
\end{equation}
Similarly to Section~3,
it suffices to consider the case where $T=1$ and $D=1$.
Assume \eqref{eq:5.3}.
Then $(U,V)$ is defined for almost all $x\in{\bf R}^N$,
$t\in(0,t_*)$ and $\tau\in(0,\rho^2]$, where $t_*$ is as in \eqref{eq:5.4}.
Furthermore, $(U,V)$ satisfies \eqref{eq:3.5} and \eqref{eq:3.7}
for almost all $x\in{\bf R}^N$, $t\in(0,t_*)$ and $\tau\in(0,\rho^2]$.

Set $W(x,t) := U(x,t) + V(x,t)$ and $M_{z,\tau}(w):=M_{z,\tau}(u) + M_{z,\tau}(v)$.
Then we have:
\begin{lemma}
\label{Lemma:6.1}
Assume \eqref{eq:5.3} and \eqref{eq:6.1}.
Let $(U,V)$ satisfy \eqref{eq:3.5}.
Assume that there exist constants $a>0$, $b\ge 1$ and $c\ge 0$ such that
\begin{equation}
\label{eq:6.3}
W(x,t)\ge a G\left( x, \frac{t+\rho^2}{b}\right)
\left[ \log \frac{t}{\rho^2} \right]^c
\end{equation}
for almost all $x\in \mathbf{R}^N$ and $t\in[\rho^2,t_*)$.
Then there exists a positive constant $\gamma_3$ depending only on $N$, $D_1$, $D_2$, $p$ and $q$ such that
$$
W(x,t)\ge
\gamma_3\frac{a^p b^{\frac{N}{2}(p-2)}}{pc+1}
G \left( x, \frac{t+\rho^2}{pb}\right)
\left[ \log \frac{t}{\rho^2} \right]^{pc+1}
$$
for almost all $x\in \mathbf{R}^N$ and $t\in[\rho^2,t_*)$.
\end{lemma}
\textbf{Proof.}
Since $p>1$, it follows that
$W(x,t)^p\le C[U(x,t)^p + V(x,t)^p]$ for almost all $x\in{\bf R}^N$ and $t\in(0,t_*)$.
Then, by \eqref{eq:3.5} and \eqref{eq:6.3} we have
\begin{equation*}
\begin{split}
W(x,t) & \ge (D')^{-\frac{N}{2}}\int_0^t\int_{{\bf R}^N}G(x-y,t-s)[U(y,s)^p+V(y,s)^p]\,dy\,ds\\
 & \ge C\int_{\rho^2}^t\int_{{\bf R}^N}G(x-y,t-s)W(y,s)^p\,dy\,ds\\
 & \ge Ca^p\int_{\rho^2}^t\int_{{\bf R}^N}
    \left[\log\frac{s}{\rho^2}\right]^{pc}G(x-y,t-s) G\left(y,\frac{s+\rho^2}{b}\right)^p\,dy\,ds
\end{split}
\end{equation*}
for almost all $x\in{\bf R}^N$ and $t\in[\rho^2,t_*)$.
This together with Lemma~\ref{Lemma:2.2} and \eqref{eq:6.2} implies that
\begin{equation*}
\begin{split}
W(x,t) & \ge C a^pb^{\frac{N}{2}(p-2)}G\left( x,\frac{t+\rho^2}{pb}\right)
\int_{\rho^2}^t (s+\rho^2)^{-1}\left[\log\frac{s}{\rho^2}\right]^{pc}\,ds\\
  & \ge Ca^p b^{\frac{N}{2}(p-2)}(pc+1)^{-1}
  G\left( x,\frac{t+\rho^2}{pb}\right)\left[ \log\frac{t}{\rho^2} \right]^{pc+1}
\end{split}
\end{equation*}
for almost all $x\in \mathbf{R}^N$ and $t\in[\rho^2,t_*)$.
Thus Lemma~\ref{Lemma:6.1} follows.
$\Box$
\vspace{5pt}
\newline
\textbf{Proof of assertion~(2) of Theorem~\ref{Theorem:1.3}.}
It suffices to prove assertion~(2) under the assumption that $T=1$ and $D=1$.
Let $z\in{\bf R}^N$, $\rho\in(0,1/\sqrt{10})$ and $\tau\in(0,\rho^2]$.
By \eqref{eq:3.7} we have
\begin{equation*}
  W(x,t) \ge c_* M_{z,\tau}(w) G(x,t+\rho^2)
\end{equation*}
for almost all $x\in \mathbf{R}^N$, $t\in[0,t_*)$ and $\tau\in(0,\rho^2]$.
Then, by Lemma~\ref{Lemma:6.1} we obtain
\begin{equation}
\label{eq:6.4}
\infty>W(x,t)\ge a_n G\left( x, \frac{t+\rho^2}{b_n}\right)\left[ \log \frac{t}{\rho^2} \right]^{c_n},
\qquad n=0,1,2,\dots,
\end{equation}
for almost all $x\in \mathbf{R}^N$, $t\in[\rho^2,t_*)$ and $\tau\in(0,\rho^2]$.
Here $\{a_n\}_{n=0}^\infty$, $\{b_n\}_{n=0}^\infty$ and $\{c_n\}_{n=0}^\infty$ are sequences defined by
\begin{equation}
\label{eq:6.5}
a_{n+1}:=\gamma_3\frac{a_n^p b_n^{\frac{N}{2}(p-2)}}{pc_n+1},
\qquad
b_{n+1}:=pb_n,
\qquad
c_{n+1}:=pc_n + 1,
\end{equation}
for $n=0,1,2,\dots$ with
$$
a_0:=c_* M_{z,\tau}(w),\qquad b_0:=1,\qquad c_0:=0,
$$
where $\gamma_3$ is as in Lemma~\ref{Lemma:6.1}.
Then we have
\begin{equation}
\label{eq:6.6}
b_n = p^n,
\qquad
c_n = \frac{p^n-1}{p-1},
\end{equation}
for $n=0,1,2,\dots$.
Furthermore, similarly to \eqref{eq:3.12},
by \eqref{eq:6.5} we can find $a_*>0$ such that
\begin{equation}
\label{eq:6.7}
a_n\ge a_*^{p^n}M_{z,\tau} (w)^{p^n}
\end{equation}
for $n=0,1,2,\dots$.
By \eqref{eq:6.4}, \eqref{eq:6.6} and \eqref{eq:6.7}
we see that
\begin{equation*}
\begin{split}
\infty>W(x,t) & \ge a_*^{p^n} M_{z,\tau}(w)^{p^n} G\left(x,\frac{t+\rho^2}{p^n}\right)
\left[\log\frac{t}{\rho^2} \right]^{\frac{p^n-1}{p-1}}\\
  & =\left\{ a_* M_{z,\tau}(w)
    \left[\log\frac{t}{\rho^2} \right]^\frac{1}{p-1} \right\}^{p^n}
    G\left(x,\frac{t+\rho^2}{p^n}\right)
    \left[ \log \frac{t}{\rho^2} \right]^{-\frac{1}{p-1}}
  \end{split}
\end{equation*}
for almost all $x\in \mathbf{R}^N$, $t\in[\rho^2,t_*)$ and $\tau\in(0,\rho^2]$ and for all $n=0,1,2,\dots$.
This implies that
$$
a_* M_{z,\tau}(w)\left[\log\frac{t}{\rho^2} \right]^\frac{1}{p-1} \le 1
$$
for almost all $t\in[\rho^2,t_*)$ and $\tau\in(0,\rho^2]$.
Since $t_*>1/2$ (see \eqref{eq:5.4}), we deduce that
$$
M_{z,\tau}(u) + M_{z,\tau}(v) \le
a_*^{-1} \left[ \log\frac{1}{2\rho^2} \right]^{-\frac{N}{2}}
\le C\left[\log\biggr(e+\frac{1}{\rho}\biggr)\right]^{-\frac{N}{2}}
$$
for all $z\in \mathbf{R}^N$ and $\rho\in(0,1/\sqrt{10})$ and for almost all $\tau\in(0,\rho^2]$.
Then, applying the same argument as in the proof of Theorem~\ref{Theorem:1.1},
we obtain the desired result, and the proof of assertion~(2) of Theorem~\ref{Theorem:1.3} is complete.
$\Box$
\section{Initial traces in cases (D) and (E)}
In this section we focus on cases (D) and (E), that is,
$$
0<p<q,\qquad  pq>1,\qquad
\frac{q+1}{pq-1}>\frac{N}{2},
\qquad q\ge 1+\frac{2}{N},
$$
and prove assertions~(3) and (4) of Theorem~\ref{Theorem:1.3}.
\vspace{5pt}
\newline
\textbf{Proof of assertions~(3) and (4) of Theorem~\ref{Theorem:1.3}.}
Similarly to Section~3, it suffices to consider the case where $T=1$ and $D=1$.
Furthermore, by Theorem~\ref{Theorem:1.1} (see \eqref{eq:1.9})
we have only to prove
\begin{equation}
\label{eq:7.1}
\sup_{z\in \mathbf{R}^N}\int_0^1 \tau^{-1}\left[ \frac{\mu(B(z,\tau))}{\tau^{N-\frac{N+2}{q}}} \right]^q\, d\tau\le\gamma,
\end{equation}
where $\gamma$ is a positive constant depending only on
$N$, $D_1$, $D_2$, $p$ and $q$.

Let $\rho_*=1/5$.  By \eqref{eq:3.17} with $\tau=\rho^2$ we have
$$
\sup_{z\in{\bf R}^N}\int_{B(z,\rho)} v(x,\rho^2)\,dx \le\gamma'
$$
for almost all $\rho\in(\rho_*/2,\rho_*)$.
Here $\gamma'$ is a positive constant depending only on $N$, $D_1$, $D_2$, $p$ and $q$.
Then Lemma~\ref{Lemma:2.4} implies that
\begin{equation}
\label{eq:7.2}
\sup_{z\in{\bf R}^N}\int_{B(z,\sqrt{2}\rho)} v(x,\rho^2)\,dx \le C\gamma'
\end{equation}
for almost all $\rho\in(\rho_*/2,\rho_*)$.
Furthermore,
similarly to \eqref{eq:5.13} and \eqref{eq:5.14}, we have
\begin{equation}
\label{eq:7.3}
\begin{split}
 &  \int_{B(z,\sqrt{2}\rho)} v(x,\rho^2) \, dx\\
 & \ge C\int_0^{\rho^2}\int_{B(z,\sqrt{s})}\int_{B(z,\sqrt{2}\rho)}
   G(x-y, \rho^2-s) [S(s)\mu](y)^q \,dx\,dy\,ds\\
 & \ge C\int_0^{\rho^2} \int_{B(z,\sqrt{s})} [S(s)\mu](y)^q \, dy \, ds
\end{split}
\end{equation}
for all $z\in{\bf R}^N$ and almost all $\rho\in(\rho_*/2,\rho_*)$.
Combining \eqref{eq:7.2} and \eqref{eq:7.3} and applying a similar argument as in \eqref{eq:5.15} and \eqref{eq:5.16},
we obtain
\begin{align*}
  C\gamma' & \ge \sup_{z\in{\bf R}^N}\int_0^{\rho^2} \int_{B(z,\sqrt{s})} [S(s)\mu](y)^q \, dy \, ds\\
   & \ge C\sup_{z\in{\bf R}^N}\int_0^{\rho^2}\mu(B(z,\sqrt{s}))^q
  \int_{B(z,\sqrt{s})} G\left(y-z,\frac{s}{2}\right)^q\,dy\,ds\\
   & \ge C\sup_{z\in{\bf R}^N}\int_0^{\rho^2} \mu(B(z,\tau))^q \tau^{-Nq+N+1}\, d\tau
   =C\sup_{z\in{\bf R}^N}\int_0^{\rho^2}\tau^{-1}\left[\frac{\mu(B(z,\tau))}{\tau^{N-\frac{N+2}{q}}}\right]^q \, d\tau
\end{align*}
for all $\rho\in(\rho_*/2,\rho_*)$.
This implies that
$$
\sup_{z\in{\bf R}^N}\int_0^{\rho_*^2}\tau^{-1}\left[\frac{\mu(B(z,\tau))}{\tau^{N-\frac{N+2}{q}}}\right]^q \, d\tau\le C.
$$
Then, applying a similar argument as in \eqref{eq:5.20}, we obtain \eqref{eq:7.1}.
Thus assertions~(3) and (4) of Theorem~\ref{Theorem:1.3} follows.
Therefore the proof of Theorem~\ref{Theorem:1.3} is complete.
$\Box$
\section{Proofs of Corollaries~\ref{Corollary:1.1} and \ref{Corollary:1.2}}
{\bf Proof of Corollary~\ref{Corollary:1.1}.}
The proof is by contradiction.
Assume that problem~\eqref{eq:1.1} possesses a global-in-time nontrivial solution.
Then, by Definition~\ref{Definition:1.1}
we see that problem~\eqref{eq:1.1} possesses a global-in-time positive solution~$(u,v)$.

Let $\tau>0$ and set $\tilde{u}(x,t):=u(x,t+\tau)$ and $\tilde{v}(x,t):={v}(v,t+\tau)$.
Then, by Definition~\ref{Definition:1.1}, for almost all $\tau>0$,
we see that $(\tilde{u},\tilde{v})$ is a global-in-time positive solution to \eqref{eq:1.1}
with $(\tilde{u}(0),\tilde{v}(0))=(u(\tau),v(\tau))$.
Furthermore, we observe from Theorem~\ref{Theorem:1.2} that
$(u(\tau),v(\tau))$ is the initial trace of the solution~$(\tilde{u},\tilde{v})$.
Therefore, by Theorem~\ref{Theorem:1.1} we can find $\gamma>0$ such that
\begin{equation}
\label{eq:8.1}
\int_{B(0,\sigma)}\tilde{u}(x,0)\,dx\le\gamma\sigma^{N-\frac{2(p+1)}{pq-1}},
\quad
\int_{B(0,\sigma)}\tilde{v}(x,0)\,dx\le\gamma\sigma^{N-\frac{2(q+1)}{pq-1}},
\end{equation}
for all $\sigma>0$.

Assume that
\begin{equation}
\label{eq:8.2}
\frac{q+1}{pq-1}>\frac{N}{2}.
\end{equation}
Then it follows from \eqref{eq:8.1} that
$$
\lim_{\sigma\to\infty}\int_{B(0,\sigma)}v(x,\tau)\,dx=\lim_{\sigma\to\infty}\int_{B(0,\sigma)}\tilde{v}(x,0)\,dx=0.
$$
Since $\tau$ is arbitrary, we see that $v(x,t)=0$ for almost all $x\in{\bf R}^N$ and $t\in(0,\infty)$.
This contradicts that $(u,v)$ is a positive solution.
Thus problem~\eqref{eq:1.1} possesses no global-in-time nontrivial solutions under assumption~\eqref{eq:8.2}.

It remains to consider the case of
\begin{equation}
\label{eq:8.3}
\frac{q+1}{pq-1}=\frac{N}{2}.
\end{equation}
We consider case (B).
Then, by Theorem~\ref{Theorem:1.3} we can find $\gamma'>0$ such that
$$
\int_{B(0,\sigma)}\tilde{v}(x,0)\,d\tau
\le\gamma'\biggr[\log\left(e+\frac{T^\frac{1}{2}}{\sigma}\right)\biggr]^{-\frac{1}{pq-1}}
$$
for all $0<\sigma<T^{\frac{1}{2}}$ and all $T>0$.
Let $T>1$ and set $\sigma=T^{\frac{1}{4}}$. Then we see that
$$
\int_{B(0,T^{\frac{1}{4}})}v(x,\tau)\,dx=
\int_{B(0,T^{\frac{1}{4}})}\tilde{v}(x,0)\,dx
\le\gamma'\biggr[\log\left(e+T^{\frac{1}{4}}\right)\biggr]^{-\frac{1}{pq-1}}\to 0
$$
as $T\to\infty$.
Similarly to the case of \eqref{eq:8.2},
we deduce that $v(x,t)=0$ for almost all $x\in{\bf R}^N$ and $t\in(0,\infty)$,
which is a contradiction.
Case~(C) can be treated in the same manner.
Therefore problem~\eqref{eq:1.1} possesses no global-in-time nontrivial solutions under assumption~\eqref{eq:8.3}.
Thus Corollary~\ref{Corollary:1.1} follows.
$\Box$\vspace{5pt}
\newline
{\bf Proof of Corollary~\ref{Corollary:1.2}.}
The proof is by contradiction.
Assume that problem~\eqref{eq:1.1} with \eqref{eq:1.2} possesses a solution~$(u,v)$
in ${\bf R}^N\times(0,T)$ for some $T>0$.
Theorem~\ref{Theorem:1.2} implies that the initial trace of $(u,v)$ coincides with the initial data of $(u,v)$.
\newline
Case (A): It follows from \eqref{eq:1.11} that
$$
\int_{B(0,\sigma)}\mu(x)\,dx\ge Cc_{a,1}\sigma^{N-\frac{2(p+1)}{pq-1}},
\quad
\int_{B(0,\sigma)}\nu(x)\,dx\ge Cc_{a,2}\sigma^{N-\frac{2(q+1)}{pq-1}},
$$
for all sufficiently small $\sigma>0$.
This implies that \eqref{eq:1.9} does not hold if either $c_{a,1}$ or $c_{a,2}$ is sufficiently large.
This contradicts Theorem~\ref{Theorem:1.1}. Thus assertion~(a) follows.
\newline
Case (B):
It follows from \eqref{eq:1.12} that
$$
\int_{B(0,\sigma)}\mu(x)\,dx\ge Cc_{b,1}\sigma^{N-\frac{2(p+1)}{pq-1}}|\log\sigma|^{-\frac{p}{pq-1}},
\quad
\int_{B(0,\sigma)}\nu(x)\,dx\ge Cc_{b,2}|\log\sigma|^{-\frac{1}{pq-1}},
$$
for all sufficiently small $\sigma>0$.
Then
$$
\int_0^\sigma \tau^{-1}\left[\frac{\mu(B(0,\tau))}{\tau^{N-\frac{2(p+1)}{pq-1}}}\right]^q\,d\tau
\ge C^qc_{b,1}^q\int_0^\sigma \tau^{-1}|\log\tau|^{-1-\frac{1}{pq-1}}\,d\tau
\ge Cc_{b,1}^q|\log\sigma|^{-\frac{1}{pq-1}}
$$
for all sufficiently small $\sigma>0$.
Then these estimates contradict assertion~(1) of Theorem~\ref{Theorem:1.3}
if either $c_{b,1}$ or $c_{b,2}$ is sufficiently large.
Thus assertion~(b) follows.
\newline
Case (C): It follows from \eqref{eq:1.13} that
$$
\int_{B(0,\sigma)}\mu(x)\,dx\ge Cc_{c,1}|\log\sigma|^{-\frac{N}{2}},
\quad
\int_{B(0,\sigma)}\nu(x)\,dx\ge Cc_{c,2}|\log\sigma|^{-\frac{N}{2}},
$$
for all sufficiently small $\sigma>0$.
Then these estimates contradict assertion~(2) of Theorem~\ref{Theorem:1.3}
if either $c_{c,1}$ or $c_{c,2}$ is sufficiently large.
Thus assertion~(c) follows.
\newline
Case (D):
Let $\epsilon>0$ be as in assertion~(d).
It follows from $q>1+2/N$ that $(N+2)/q<N$.
Then, by \eqref{eq:1.14} we have
\begin{equation*}
\begin{split}
\int_{B(0,\sigma)}\mu(x)\,dx
 & \ge C\int_0^\sigma r^{-\frac{N+2}{q}+N-1+\epsilon}r^{-\epsilon} h_1(r)\,dr\\
 & \ge C\sigma^{-\epsilon} h_1(\sigma)\int_0^\sigma r^{-\frac{N+2}{q}+N-1+\epsilon}\,dr
\ge C\sigma^{N-\frac{N+2}{q}}h_1(\sigma)
\end{split}
\end{equation*}
for all sufficiently small $\sigma>0$. This implies that
\begin{equation}
\label{eq:8.4}
\int_0^t
\tau^{-1}\left[\frac{\mu(B(x,\tau))}{\tau^{N-\frac{N+2}{q}}}\right]^q\,d\tau
\ge\int_0^t \tau^{-1}h_1(\tau)^q\,d\tau
\end{equation}
for all sufficiently small $t>0$.
On the other hand,
since $h_1$ is positive continuous function in $(0,1)$,
if
$$
\int_0^1 \tau^{-1}h_1(\tau)^q\,d\tau=\infty,
$$
then
\begin{equation}
\label{eq:8.5}
\int_0^t \tau^{-1}h_1(\tau)^q\,d\tau=\infty\quad\mbox{for $t\in(0,1)$}.
\end{equation}
Furthermore, for any $r>0$,
by Lemma~\ref{Lemma:2.4} we have
\begin{equation}
\label{eq:8.6}
\sup_{x\in{\bf R}^N}\nu(B(x,1))\le C\sup_{x\in{\bf R}^N}\nu(B(x,r)).
\end{equation}
Then, under the assumptions of assertion~(d), we see a contradiction
between \eqref{eq:8.4}, \eqref{eq:8.5}, \eqref{eq:8.6} and assertion~(3) of Theorem~\ref{Theorem:1.3}.
Thus assertion~(d) follows.
\newline
Case~(E): Since $q=1+2/N$, by \eqref{eq:1.16} we have
$$
\int_{B(0,\sigma)}\mu(x)\,dx\ge C\int_0^\sigma \tau^{-1}h_2(\tau)\,d\tau
$$
for all sufficiently small $\sigma>0$.
This implies that
\begin{equation}
\label{eq:8.7}
\int_0^t \tau^{-1}\mu(B(x,\tau))^q\,d\tau
\ge C\int_0^t r^{-1}\left[\int_0^r \tau^{-1}h_2(\tau)\,d\tau\right]^q\,dr
\end{equation}
for all sufficiently small $t>0$.
On the other hand,
since $h_2$ is positive continuous function in $(0,1)$, if
$$
\int_0^1\left[\int_0^r \tau^{-1}h_2(\tau)\,d\tau\right]^qr^{-1}\,dr=\infty,
$$
then
\begin{equation}
\label{eq:8.8}
\int_0^t\left[\int_0^r \tau^{-1}h_2(\tau)\,d\tau\right]^qr^{-1}\,dr=\infty
\quad\mbox{for $t\in(0,1)$}.
\end{equation}
Furthermore, \eqref{eq:8.6} holds.
Therefore, under the assumptions of assertion~(e), we see a contradiction
between \eqref{eq:8.6}, \eqref{eq:8.7}, \eqref{eq:8.8} and assertion~(4) of Theorem~\ref{Theorem:1.3}.
Thus assertion~(e) follows.
Furthermore, assertion~(f) also follows from Lemma~\ref{Lemma:2.4} and Theorem~\ref{Theorem:1.1}
(see also Remark~\ref{Remark:1.2}).
Therefore the proof of Corollary~\ref{Corollary:1.2} is complete.
$\Box$
\medskip

\noindent
{\bf Acknowledgements.}
The first author was supported partially by the Grant-in-Aid for Early-Career Scientists (No.~19K14569).
The second author of this paper was supported in part
by the Grant-in-Aid for Scientific Research (S)(No.~19H05599)
from Japan Society for the Promotion of Science.
\bibliographystyle{amsplain}


\end{document}